\documentclass{amsart}
\usepackage{latexsym}
\usepackage{graphicx}
\usepackage{amscd}
\usepackage{pstricks,pst-node}
\newfont{\msam}{msam10}
\newcommand{\la}{\label}
\newtheorem{lemma}{Lemma}[section]
\newtheorem{proposition}[lemma]{Proposition}
\newtheorem{theorem}[lemma]{Theorem}
\newtheorem{corollary}[lemma]{Corollary}

\theoremstyle{definition}

\newtheorem{definition}[lemma]{Definition}

{

}
\theoremstyle{remark}

\def\bs#1{\boldsymbol{#1}}
\def\ms#1{\mathcal{#1}}

\def\Ta{\mbox{\sf Tails}(\bs{A})}

\def\ta{\mbox{\sf tails}(\bs{A})}

\def\Gma{\mbox{\sf GrMod}(\bs{A})}
\def\gma{\mbox{\sf grmod}(\bs{A})}
\def\pr#1{\mbox{\sf proj}(#1)}

\def\P{\mathbb{P}_{q}^{2}}
\def\Z{\mathbb{Z}}
\def\U{\mathfrak{U}}
\def\H{\underline{\check{H}}{}}
\def\X{\mathbb{X}}
\def\Y{\mathbb{Y}}
\def\I{\mathbb{I}}
\def\vv{\mathfrak{v}}
\def\u{\mathfrak{u}}
\def\c{\mathbb{C}}
\def\PP{\mathbb{P}^{1}}
\def\Qc#1{\mbox{\sf Qcoh}(#1)}
\def\C#1{\mbox{\sf coh}(#1)}

\def\Hom\sub#1#2{\underline{\mbox{\rm Hom}}_{\bs{A}}(#1,#2)}
\def\Ext\sub#1#2{\underline{\mbox{\rm Ext}}^{i}_{\bs{A}}(#1,#2)}
\numberwithin{equation}{section}
\let\cal\mathcal

\def\Cscr{{\cal C}}

\def\Escr{{\cal E}}

\def\Mscr{{\cal M}}

\def\Oscr{{\cal O}}

\def\Rscr{{\cal R}}

\def\Xscr{{\cal X}}
\def\Yscr{{\cal Y}}

\let\blb\mathbb
\def\CC{{\blb C}} 
\def\XX{{\blb X}} 
\def\YY{{\blb Y}}

\def \PP{{\blb P}}

\def \ZZ{{\blb Z}}

\def\opp{\operatorname{opp}}

\newcommand{\proj}{\operatorname{proj}}

\def\Id{\operatorname{Id}}

\def\Lotimes{\overset{\textbf{L}}{\otimes}}

\def\Mod{\operatorname{Mod}}
\def\mod{\operatorname{mod}}

\def\coh{\mathop{\text{\upshape{coh}}}}

\def\Ext{\operatorname {Ext}}
\def\Hom{\operatorname {Hom}}
\def\End{\operatorname {End}}
\def\RHom{\operatorname {\textbf{R}Hom}}

\def\Tor{\operatorname {Tor}}
\def\End{\operatorname {End}}

\def\rk{\operatorname {rk}}

\def\r{\rightarrow}

\DeclareMathOperator{\Aut}{Aut}

\newdimen\uboxsep \uboxsep=1ex
\def\uboxn#1{\vtop to 0pt{\hrule height 0pt depth 0pt\vskip\uboxsep
\hbox to 0pt{\hss #1\hss}\vss}}

\def\uboxs#1{\vbox to 0pt{\vss\hbox to 0pt{\hss #1\hss}
\vskip\uboxsep\hrule height 0pt depth 0pt}}

\begin{document}
\title[Ideal Classes of the Weyl Algebra]
{Ideal Classes of the Weyl Algebra\\ and 
Noncommutative Projective Geometry\\ 
(with an Appendix by Michel Van den Bergh)}
%
%
	\author{Yuri Berest}
\address{Department of Mathematics,
Cornell University, Ithaca, NY 14853-4201, USA}
\email{berest@math.cornell.edu}
%
%
\thanks{Research supported in part by NSF grant DMS 00-71792.}
%
%
     \author{George Wilson}
     \address{Department of Mathematics, 
      Imperial College, London SW7 2BZ, UK}
     \email{g.wilson@ic.ac.uk}
\begin{abstract}
Let $\, \mathfrak{R} \,$ be the set of isomorphism classes
of right ideals in the Weyl algebra $\, A = A_{1}(\c) \,$, and let 
$\, \mathfrak{C} \,$ be the set of isomorphism classes of
triples $\, (V;\, \X,\, \Y) \,$, where $\, V \,$ is a 
finite-dimensional (complex) vector space, and 
$\, \X,\, \Y \,$ are endomorphisms of $\, V \,$ such that
$\, [\X,\Y] + \I \,$ has rank $ 1\,$. Following a suggestion
of L.~Le Bruyn, we define a map $\, \theta:\, \mathfrak{R} 
\to \mathfrak{C} \,$ by appropriately extending an 
ideal of $\, A \,$ to a sheaf over a quantum projective 
plane, and then using standard methods of homological algebra.
We prove that $\, \theta\,$ is inverse to a bijection 
$\, \omega:\, \mathfrak{C} \to \mathfrak{R} \,$ constructed
in \cite{BW} by a completely different method. The main step 
in the proof is to show that $\,\theta \,$ is equivariant
with respect to natural actions of the group $\, G = 
\mbox{\rm Aut}(A)\,$ on 
$\,\mathfrak{R} \,$ and $\,\mathfrak{C} \,$: for that we have 
to study also the extensions of an ideal to certain {\it 
weighted} quantum projective planes. Along the way, we find
an elementary description of $\, \theta\,$.
\end{abstract}
\maketitle

\section{Introduction}
\la{Intro}

This is a sequel to our earlier paper \cite{BW}; however, it can, and
probably should, be read independently of that work. 
We first recall the main results of \cite{BW}.

Let $ A $ be the complex Weyl algebra; that is, $ A $ is the associative 
algebra over $ \mathbb{C} $ generated by two elements $ x $ and $ y $ subject 
to the relation $\, [x,y]=1\,$. Let $ \mathfrak{R} $ be the set
of isomorphism classes of finitely generated torsion-free rank one
right $A$-modules. Since $ A $ is Noetherian, each right ideal of 
$ A $ is a module of this kind; conversely, because $ A $ has a 
quotient (skew) field, it is easy to see that each such 
module $ M $ is isomorphic to a right ideal of $ A \,$. For short,
we shall often refer to a module $ M $ as an ``ideal'', even when
we do not have in mind any particular embedding of $ M $ in $ A\,$.
Let $ G $ be the automorphism group of $ A\,$; there is a natural
action of $ G $ on $ \mathfrak{R} \,$. Finally, for each 
$\, n \geq 0\,$, let $ \mathfrak{C}_n $ be the space of equivalence 
classes (modulo simultaneous conjugation) of pairs $\, (\X,\Y) \,$
of  $ n \times n $ matrices over $ \mathbb{C} $ such 
that\footnote{In \cite{BW} we worked with the space of pairs such
that $\, [\X,\Y] - \I \ \mbox{has rank}\ 1 \,$; here we
identify these spaces via the map $\, (\X,\Y) \leftrightarrow 
(\X^{t},\Y^{t}) \,$.  When $\, n = 0\,$, the space $ \mathfrak{C}_n $ 
is supposed to be a point: 
as in (\ref{I1}), we shall sometimes disregard this case.}
\begin{equation}
\label{I1}
[\X,\Y] + \I \quad \mbox{has rank}\ 1 \ .
\end{equation}
For brevity, we shall often refer to a point of  $ \mathfrak{C}_n $ 
simply as a ``pair of matrices''. There is a natural action of 
$ G $ on each space $ \mathfrak{C}_n \,$; it is obtained, roughly
speaking, by thinking of the pairs $\, (\X,\Y) \,$ as points dual
to the coordinate functions $ x $ and $ y $ that generate $ A\,$
(see Section~\ref{Sect7} below for the precise definition).
In \cite{BW} we showed that this action is transitive. 
Let  $ \mathfrak{C} $ be the (disjoint) union of the spaces
$ \mathfrak{C}_n \,$. The main result of \cite{BW} was the 
following
\begin{theorem}
\la{IT1}
There is a bijective map $\, \omega: \mathfrak{C} \to \mathfrak{R} \,$
which is equivariant with respect to the action of $ G\,$.
\end{theorem}

Part of the significance of this Theorem becomes clear if we think
of the Weyl algebra as a noncommutative deformation of the
polynomial ring $ \mathbb{C}[x,y]\,$. In this case the 
analogue of Theorem~\ref{IT1} is elementary, because each
isomorphism class of ideals has a unique representative
as an ideal of finite codimension; thus in the commutative 
case, the analogue of our space $ \mathfrak{R} $ is the 
disjoint union of the point Hilbert schemes $\,\mbox{Hilb}_n(\mathbb{A}^2)\,$
of the affine plane. As is well known (and almost tautological),
$\,\mbox{Hilb}_n(\mathbb{A}^2)\,$ can be identified with the space
of (equivalence classes of) pairs $\, (\X,\Y)\,$ of 
{\it commuting} $\, n \times n \,$  matrices possessing a cyclic vector:
to an ideal $\, I\,$ of finite codimension we assign the pair 
$\, (\X,\Y)\,$ of maps on the quotient $\, \c[x,y]/I \,$ induced by 
multiplication by $\, x \,$ and $\, y \,$. According to Nakajima 
(see \cite{N}), our space 
$ \mathfrak{C}_n $ can be obtained from 
$\,\mbox{Hilb}_n(\mathbb{A}^2)\,$ by a deformation
of complex structure. In commutative algebraic geometry it is a 
basic principle that a (flat) deformation of varieties should give 
rise to a deformation of any reasonable moduli space of bundles
(or coherent sheaves); Theorem~\ref{IT1} suggests that this principle
should extend also to noncommutative deformations. However, at the 
present time even the expression ``moduli space'' seems not to have any
precise meaning in the noncommutative case; our space 
$ \mathfrak{R} \,$, for example, does not to our knowledge possess
any {\it intrinsic} algebraic structure, which is why we referred
to it above simply as a set.

The description of the map $\, \omega: \mathfrak{C} \to \mathfrak{R} \,$
given in \cite{BW} is not direct, but passes through a third space,
the {\it adelic Grassmannian} $ \mbox{Gr}^{\mbox{\scriptsize \rm ad}} $
that parametrizes rational solutions of a certain integrable system
(the KP hierarchy). Indeed, in \cite{BW} we defined $ \omega $
to be the composition of a bijection $\, \beta: 
\mathfrak{C} \to  \mbox{Gr}^{\mbox{\scriptsize \rm ad}} \,$
constructed in \cite{W1} and a bijection
$\, \alpha: \mbox{Gr}^{\mbox{\scriptsize \rm ad}} \to
\mathfrak{R} \,$ constructed by Cannings and Holland
in \cite{CH}. Theorem~\ref{IT1} was then proved by following
through what happens to the natural action of $ G $ on
$ \mathfrak{R} $ under the bijections $ \alpha $ and 
$ \beta \,$, a tricky process, since the action 
of $ G $ on $ \mbox{Gr}^{\mbox{\scriptsize \rm ad}} $ is difficult to
describe. In any case, although $ \mbox{Gr}^{\mbox{\scriptsize \rm ad}} $
is an interesting object in its own right, it is hard not to feel
that it is {\it de trop} in the question of classifying ideals
of $ A\,$. Another imperfection in \cite{BW} is that there we gave 
no explicit description of the inverse map to 
$ \omega \,$ (Cannings and Holland do indeed explain what is the
inverse of their bijection $ \alpha \,$, but no description was known 
for the inverse of the map $ \beta \,$). For these (and other) reasons,
we wish to rederive Theorem~\ref{IT1} in a way that makes no reference
to $ \mbox{Gr}^{\mbox{\scriptsize \rm ad}} \,$.

To that end, we take up an idea of L.~Le Bruyn \cite{LeB}. If we think
of $ A $ as the ring of functions on a ``quantum affine plane''
$ \mathbb{A}_{q}^{2} \,$, then an $ A$-module $ M $ (in particular, an
ideal) is to be thought of as a vector bundle 
(or coherent sheaf) over $ \mathbb{A}_{q}^{2} \,$. Le Bruyn's
idea is to extend $ M $ to a sheaf $ \ms{M} $ over the quantum
{\it projective} plane $ \P \,$, and then use a noncommutative
version of Barth's classification \cite{Bar} of bundles over
$ \mathbb{P}^2 $ to obtain algebraic invariants of $ M\,$.
Here $ \P $ is taken in the sense of M.~Artin \cite{A}:
following Serre's classic paper \cite{Se}, Artin suggests to
define a (noncommutative) projective variety via its
homogeneous coordinate ring $ \bs{A} \,$, so that a sheaf over
such a variety is represented by a graded $ \bs{A}$-module
(modulo finite-dimensional modules). We give a quick sketch of this theory
in Section~\ref{Sect2} below: readers who are not familiar with
it may understand the rest of this Introduction by analogy with
the commutative case. The homogeneous coordinate ring of $ \P $
is the graded ring of ``noncommutative polynomials''
$\, \bs{A} = \mathbb{C}[X,Y,Z] \,$, where the generators
$ X, Y $ and $ Z $ all have degree $ 1\,$, $\, Z \,$ commutes
with $ X $ and $ Y $, and $\, [X,Y]= Z^2\,$. As for the 
classification of bundles over $ \P \,$, Le Bruyn suggests
to use the version of Beilinson \cite{B}, which goes very smoothly
in the noncommutative case (cf. \cite{Bon}).
The following consequence of that theory is all that will 
concern us in this paper. Suppose we have a sheaf $ \ms{F} $
over $ \P $ satisfying the vanishing conditions
\begin{equation}
\label{I2}
H^{i}(\P, \ms{F}(-2)) = H^{i}(\P, \ms{F}(-1)) = 
H^{i}(\P, \ms{F}) = 0 \quad \mbox{for all}\ i \not= 1\ .
\end{equation}
Then $ \ms{F} $ is determined by the representation
\begin{equation}
\label{I3}
H^{1}(\P, \ms{F}(-2))\  
\begin{array}{c}
\longrightarrow \\*[-1ex] 
\longrightarrow \\*[-1ex]
\longrightarrow 
\end{array} 
\ H^{1}(\P, \ms{F}(-1))\  
\begin{array}{c}
\longrightarrow \\*[-1ex] 
\longrightarrow \\*[-1ex]
\longrightarrow 
\end{array}
\  H^{1}(\P, \ms{F}) \ 
\end{equation}
of the indicated quiver with $3$ vertices and $6$ arrows
(and relations reflecting the commutation relations of 
the algebra $ \bs{A} $). In (\ref{I3})
it is understood that each set of $3$ arrows is given
by multiplication by the generators $ X\,$, $ Y $ and 
$ Z $ of $\, \bs{A} \,$. It turns out that 
any\footnote{We disregard the case of the free 
$ A$-module of rank $1\,$, which corresponds to 
$ n=0 $ below.} ideal of $ A $ has
extensions $ \ms{M} $ that satisfy (\ref{I2}); in particular,
we shall see that the (unique) extension whose restriction
to the line at infinity in $ \P $ is {\it trivial} always
satisfies (\ref{I2}). From now on,  $ \ms{M} $ will denote
this extension (it is just at this point that we part company
from L.~Le Bruyn, who chooses a different extension).
Then the three vector spaces 
$\, H^{1}(\P, \ms{M}(-2)) \,$, $\, H^{1}(\P, \ms{M}(-1)) \,$ 
and $\, H^{1}(\P, \ms{M}) \,$ have dimensions $ (n,n,n-1) $
for some $ n\geq 1\,$. Furthermore, multiplication by 
$ Z $ gives a surjection 
$\, H^{1}(\P, \ms{M}(-1)) \, \mbox{\msam \symbol{016}}\, 
H^{1}(\P, \ms{M})\, $ and an isomorphism
$\, H^{1}(\P, \ms{M}(-2)) \stackrel{\sim}{\to}
H^{1}(\P, \ms{M}(-1))\,$. If we now identify
$\, H^{1}(\P, \ms{M}(-2)) \cong H^{1}(\P, \ms{M}(-1)) 
\cong V \,$ (say) via this isomorphism, then $ X $ and 
$ Y $ become endomorphisms of $ V \,$, and it is easy to see
that they satisfy the relation (\ref{I1}).               
In this way, Le Bruyn's (modified) construction gives us a map
$ \theta : \mathfrak{R} \to \mathfrak{C}\,$.
Our aim is to prove
\begin{theorem}
\la{IT2}
The map $ \theta $ is inverse to the map $ \omega $ defined 
in \cite{BW}.
\end{theorem}

In the present paper we shall not attempt a 
direct proof of this Theorem; indeed,
the machinery involved in the definitions of $ \omega $ and
$ \theta $ is so different that this appears (at first sight)
impossible. Instead, we focus on the equivariance property
of $ \omega\,$ stated in Theorem~\ref{IT1}. Because the action
of $ G $ on each space $\, \mathfrak{C}_n \,$ is transitive, the
map $\, \omega^{-1} \,$ is uniquely determined by equivariance
once we know its effect on one point in each orbit
$\, \omega(\mathfrak{C}_n) \,$. Now, there is (at least) one point
$\, M \,$ in each orbit for which it is possible to check directly that
$\, \omega^{-1}(M) = \theta(M)\,$; granting that, 
Theorem~\ref{IT2} is equivalent to
\begin{theorem}
\la{IT3}
The map $\, \theta \,$ described above is $G$-equivariant.
\end{theorem}

One might think at first that Theorem~\ref{IT3} should follow easily
from simple considerations of functoriality: however,
the difficulty arises that an automorphism of the affine
plane does not (in general) extend to a regular 
automorphism of the projective plane, but only to a 
birational automorphism. In algebraic terms, this means
that an automorphism of the Weyl algebra does not naturally
induce any map on the graded ring $ \bs{A}\,$, which is
the only object we have to work with in the noncommutative 
case. Our idea for dealing with this problem rests on
a theorem of Dixmier (see \cite{Di}), which states that the 
group $ G $ is generated by the special automorphisms
$ \Psi_{r,\lambda} $ and $ \Phi_{s,\mu} $ defined by
\begin{equation}
\la{auto}
\begin{array}{lll}
&&\Psi_{r,\lambda}(x) = x \ ,  \quad \Psi_{r,\lambda}(y) = 
y + \lambda x^r \ ;\\*[0.6ex]
&&\Phi_{s,\mu}(x) = x + \mu y^s \ ,  \quad  \Phi_{s,\mu}(y) = y\ .
\end{array}
\end{equation}
Clearly, it is enough to prove that $ \theta $ commutes with
the action of these generators. We observe that 
$ \Psi_{r,\lambda} $ and $ \Phi_{s,\mu} $ will be homogeneous
if we assign to $ x $ and $ y $ the weights $ (1,r) $ and
$ (s,1) $ (respectively); in geometrical language, this means
that these automorphisms extend to (biregular) automorphisms
of an appropriate {\it weighted projective space}. Slightly
more generally, we shall work with the weighted projective
spaces $ \P(\bs{w}) $ for any weight vector 
$\,\bs{w} =(w_1, w_2)\,$, where $ w_1 $ and $ w_2 $ are 
positive integers (it does not simplify what follows to assume
that one of them is equal to $1$). By definition, the homogeneous
coordinate ring of $ \P(\bs{w}) $ is the graded algebra
$ \bs{A}(\bs{w}) = \mathbb{C}[X,Y,Z] \,$, where $ (X,Y,Z) $ 
have degrees $ (w_1, w_2, 1)\,$, $\, Z $ commutes with $ X $ 
and $ Y \,$, and we have the commutation relation
$\,[X,Y] = Z^{|\bs{w}|}\,$ (we set $|\bs{w}| := w_1 + w_2 $).
It is not difficult to repeat all of Le Bruyn's considerations
for any weight $ \bs{w}\,$, the only difference being that the quiver
that arises is more complicated than the one we saw in (\ref{I3}).
Thus, each ideal $ M $ extends to a sheaf $ \ms{M}_{\bs{w}} $
over $ \P(\bs{w}) \,$, and we obtain a pair of matrices, say
$ (\X(\bs{w}), \Y(\bs{w})) \,$, acting on the vector space
$ V(\bs{w}) := H^{1}(\P(\bs{w}),\, \ms{M}_{\bs{w}}(-1))\,$.
If we take $ \bs{w} = (1,r) $ (or $ \bs{w} = (s,1) $), then
the automorphisms $ \Psi_{r,\lambda} $ (or $ \Phi_{s,\mu}$)
extend to automorphisms of the graded ring $ \bs{A}(\bs{w}) \,$,
so we can follow their action on the corresponding pair of 
matrices $ (\X(\bs{w}), \Y(\bs{w})) \,$ by simple functorial 
considerations. Theorem~\ref{IT3} will therefore follow 
at once from the next theorem, which is perhaps to be considered
the main result of this paper.
\begin{theorem}[Comparison theorem]
\la{IT4}
The pair of matrices $\, (\X(\bs{w}), \Y(\bs{w})) \in \mathfrak{C} \,$
corresponding to a given ideal of $ A $ is independent of the choice 
of $ \bs{w}\,$.
\end{theorem}

The proof of Theorem~\ref{IT4} that we shall present here
does not compare different weights directly, but instead compares
each pair $\, (\X(\bs{w}), \Y(\bs{w})) \,$ with yet another pair of 
matrices $ (\X,\Y) $ which we shall extract from an ideal $ M $
in an elementary way (that is, without the use of 
homological algebra). The construction imitates the elementary
treatment in the commutative case, using the representative of finite
codimension for an ideal. The Weyl algebra has no (non-trivial)
ideals of finite codimension. However, in Section~\ref{Sect5}
below we shall construct for each $\, M \,$ two 
\,{\it fractional}\, ideals $\, M_x \,$ and $\, M_y \,$ (both
isomorphic to $M$), together with embeddings $\, r_x \,$
and $\, r_y \,$ of $\, M_x \,$ and $\, M_y \,$ as linear 
subspaces of finite codimension in $\, A\,$. Of course, in that case
$\, r_x \,$ and $\, r_y \,$ cannot be $A$-module homomorphisms,
but $\, r_x \,$ will be $ \c[y]$-linear and $\, r_y \,$ will be 
$ \c[x]$-linear. On the quotient spaces $\, V_x := A/ r_x(M_x) \,$ 
and $\, V_y := A/r_y(M_y)\,$ we therefore have endomorphisms 
$\, \Y \,$ and $\, \X \,$ (respectively) induced by multiplication
by $\, y \,$ and $\, x\,$. Furthermore, we shall
construct a {\it canonical} isomorphism 
$\,\phi :\, V_x \to V_y \,$. Identifying $ V_x $ and $ V_y $
via $ \phi \,$, we thus get a pair of matrices  $ (\X,\Y) $
(as usual, defined only up to simultaneous
conjugation).
\begin{theorem}
\la{IT5}
Let $ M $ be an ideal of $ A \,$, and for each positive 
weight vector $ \bs{w} $ let $ (\X(\bs{w}), \Y(\bs{w})) $
be the corresponding pair of endomorphisms of the vector 
space $\, V(\bs{w}) = H^{1}(\P(\bs{w}),\, \ms{M}_{\bs{w}}(-1))\,$
described earlier. Then there are isomorphisms
$\, \alpha_{x}: \, V(\bs{w}) \to V_x \,$ and
$\, \alpha_{y}: \, V(\bs{w}) \to V_y \,$
taking $ \Y(\bs{w}) $ to $ \Y $ and $ \X(\bs{w}) $ to $ \X $
respectively, and making the diagram
$$
\begin{CD}
V(\bs{w}) @>\alpha_{x}>> V_x\\
@|                       @VV\phi V \\
V(\bs{w}) @>\alpha_{y}>> V_y\\
\end{CD}
$$
commutative.
\end{theorem}

In other words, the various pairs of matrices $
(\X(\bs{w}),\Y(\bs{w})) $
and  $ (\X,\Y) $ that we have assigned to an ideal $ M $ all coincide.
Theorem~\ref{IT4} is now clear, since the elementary construction
of the pair $ (\X,\Y) $ does not involve any choice of weights.

The proof of Theorem~\ref{IT5} consists of a calculation of 
$\, H^{1}(\P(\bs{w}),\, \ms{M}_{\bs{w}}(-1))\,$  using the
\v{C}ech cohomology developed in \cite{Ver}, \cite{VW1}, \cite{VW2}. 
A key point is that the quantum ``planes'' $ \P(\bs{w}) $ 
have {\it schematic dimension} (in the sense of \cite{W})
equal to $ 1 \,$, not $ 2\,$; this means (in geometrical 
language) that $ \P(\bs{w}) $ can be covered by just two
affine open sets, analogous to the $ (X,Z)$-plane and the 
$ (Y,Z)$-plane in the commutative case. 
Of course, in that case these two
affine open sets fail to catch the point $ (0:0:1) \,$;
but apparently the quantum planes do not have ``points''
to cause that kind of trouble. It turns out that our
pair $ (M_x, M_y) $ of special representatives of an ideal
is well adapted to calculating the \v{C}ech cohomology
of this $ 2$-set covering of  $ \P(\bs{w}) \,$. For 
details we refer to Section~\ref{Sect6} below.

After this article appeared as a preprint, M. Van den Bergh 
succeeded in finding direct homological proofs of our main results:
these proofs are presented in the Appendix.  In particular,
Van den Bergh proves directly that our map $\, \theta \,$ is 
bijective, whereas in the main body of the paper we see that 
only after identifying $\, \theta \,$ with $\, \omega^{-1} $: 
this proof of bijectivity thus still depends on the arguments 
from the theory of integrable systems used in \cite{BW}. It is 
interesting that methods from integrable systems and from the 
theory of derived categories appear here as alternatives to each other 
(cf.\ the question raised in the first sentence of \cite{Be}).  

The paper is organized as follows. In Section~\ref{Sect2}
we give a brief introduction to noncommutative projective 
geometry, and summarize the results we need from the 
literature on that subject. In Section~\ref{Sect3} 
we introduce the main characters in our story, the Weyl
algebra, its various ``homogenizations'' and the associated 
projective planes $\, \P(\bs{w}) \,$.
Then in Section~\ref{Sect4} we show (for arbitrary weights)
how to extract from a given ideal of $ A $ the pairs of matrices
$\, (\X(\bs{w}),\Y(\bs{w})) \in \mathfrak{C}\, $. The next section
explains the elementary construction of the pair of matrices
$ (\X,\Y) $, using the two special representatives of an ideal
of $ A\,$. Then Section~\ref{Sect6} gives the calculation of
\v{C}ech cohomology which identifies 
$ (\X(\bs{w}),\Y(\bs{w})) $ with $ (\X,\Y) \,$, and hence proves
Theorems~\ref{IT5} and~\ref{IT4}. Section~\ref{Sect7} then deduces 
our other main results, filling in some details left vague 
in the sketch above.  Section~\ref{Sect8} establishes
Beilinson's derived equivalence for the spaces 
$\, \P(\bs{w}) \,$: although, strictly speaking, this equivalence 
is not used in the main part of the paper\footnote{However, it is used in 
an essential way in the Appendix.}, it is scarcely possible 
to understand the motivation for our construction without at least 
a cursory reading of this section.  Finally, section~\ref{Sect9} 
discusses briefly the relationship of our construction
to the original one of Le Bruyn, while Section~\ref{Sect10}
explains how it fits in with the classification of bundles
over $ \P $  given in the recent paper \cite{KKO}. 
\subsection*{Acknowledgements}
We are indebted to L.~Le Bruyn and M.~Kontsevich
for advice in the early stages of the present work;
in particular, they both (independently) pointed out 
to us in October 1998 that a point of the space
$\, \mathfrak{C}_n \,$ could be identified with
a special kind of representation of the quiver
in (\ref{I3}). We thank M. Van den Bergh for kindly 
providing the Appendix. The second author is grateful
to the Mathematics Department of Cornell University
for hospitality during April 2001. 
\section{Noncommutative Projective Geometry}
\la{Sect2}

As we mentioned in the Introduction, the starting point for
noncommutative projective geometry is the following 
result of Serre (see \cite{Se}): 
the category of coherent sheaves over a (commutative) projective
variety $ X $ is equivalent to a quotient of the category of finitely
generated graded modules over the coordinate ring $ \bs{A} $ of 
$ X\,$. This latter category makes sense also for a 
noncommutative graded ring $ \bs{A} \,$.  
In this section we give a brief overview of the theory
of noncommutative projective schemes and their 
cohomology: we introduce the notation, recall definitions and 
collect some fundamental results needed for 
understanding the main part of the paper. 
Our basic reference is \cite{AZ}. 
More details on graded algebras and modules 
can be found in \cite{NV}, on abelian categories 
(including Serre quotients) in \cite{G} and \cite{GM}, 
on (Artin-Schelter) regular algebras in \cite{AS}, 
\cite{ATV} and \cite{Steph}. The ``schematic'' structure on graded
algebras and the noncommutative \v{C}ech cohomology are
introduced and discussed in \cite{Ver}, \cite{VW1}, \cite{VW2}. 
For a general overview of the 
subject we recommend the articles \cite{A}, 
\cite{Sm} and \cite{VW}.

\subsection{Graded Algebras and Modules}
\label{S21} 

We recall that an associative algebra\footnote{Throughout 
the paper we shall 
denote graded objects (algebras, modules, $\,\ldots\,$) 
by (capital) boldface letters 
($ \bs{A},\,  \bs{M},\,\, \ldots \,$) 
distiguishing them from ungraded ones 
($ A,\, M,\, \ldots $)\,.} 
$ \bs{A} $ over a field $ k $ is 
called {\it graded} (more precisely, 
$ {\mathbb Z}$-graded) if $\, \bs{A} = 
\bigoplus_{i \in  {\mathbb Z}} A_{i} \,$ 
as a $ k$-vector space and
$ A_{i} A_{j} 
\subseteq A_{i+j} $ for all  
$ i, j \in {\mathbb Z}\, $. 
We shall assume that $\, \bs{A} \,$ is  
(both left and right) Noetherian.
If $ A_{i} = 0 $ for all 
$ i < 0 $ and $ A_{0} = k $ 
we say that
$ \bs{A} $ is {\it connected}. It is 
easy to see 
that any graded connected Noetherian 
$ k$-algebra is 
{\it locally finite}, that is, 
$ \dim_{k} A_{i} < \infty $ 
for all $ i\,$.

A (right) $ \bs{A}$-module $ \bs{M} $ is
{\it graded} if it has a vector space decomposition
$ \bs{M} = \bigoplus_{i \in  {\mathbb Z}} 
M_{i} $ compatible with the $ \Z$-grading
on $\, \bs{A}\,$, that is, $\, M_{i}\, A_{j} 
\subseteq M_{i+j} \,$ for all $\, i, j \, $.
The category of all right graded modules over $ \bs{A} $
will be denoted by $ \mbox{\sf GrMod}(\bs{A}) $:
the morphisms in $ \mbox{\sf GrMod}(\bs{A}) $ 
are {\it graded} homomorphisms of degree zero. 

For each $\, n \in {\mathbb Z} \,$ we introduce two
functorial operations on graded modules 

\vspace{0.2cm}

{\it shift in grading:} 
$\ \bs{M} = \bigoplus_{i \in  {\mathbb Z}} 
M_{i} \ \mapsto \ 
\bs{M}(n) := \bigoplus_{i \in  {\mathbb Z}} 
M_{i+n}\ ; $\\*[-0.2cm]

{\it (left) truncation:} 
$\ \bs{M} = \bigoplus_{i \in  {\mathbb Z}} 
M_{i} \ \mapsto \ 
\bs{M}_{\geq n} 
:= \bigoplus_{i \geq n } 
M_{i} \ . \\*[0.2cm]
$
We say that a graded
module $ \bs{M} $ is {\it left} 
(respectively, {\it right}) {\it bounded} if 
$\, \bs{M}_{\geq n} = \bs{M} \,$ (respectively, 
$\, \bs{M}_{\geq n} = 0 $) for some 
$ n \in {\mathbb Z}\,$.

Finitely generated graded modules form a full subcategory 
in $ \mbox{\sf GrMod}(\bs{A})\,$; it is denoted by
$ \mbox{\sf grmod}(\bs{A}) $. The shift
and truncation functors on 
$ \mbox{\sf GrMod}(\bs{A}) $ preserve this 
subcategory. 
Moreover, if $ \bs{A} $ is left 
bounded (for example, connected), so are all finitely
generated graded modules over $ \bs{A} $.
Further, if $ \bs{A} $ is locally finite, then 
every object in $ \mbox{\sf grmod}(\bs{A}) $ is
locally finite as well.

A few words on homological properties of 
graded modules. First of all, 
$ \mbox{\sf GrMod}(\bs{A}) $ is an
abelian $ k$-linear category with 
enough projective and injective objects,
so for each $\, n \geq 1 \, $ we may define the functors 
$\, \mbox{Ext}^{n}_{\bs{A}}
(\bs{M}, \mbox{---})\, $ on 
$ \, \mbox{\sf GrMod}(\bs{A})\, $ as the right derived of 
$\, \mbox{Hom}_{\bs{A}}
(\bs{M}, \mbox{---}) 
\equiv \mbox{Hom}_{\mbox{\scriptsize GrMod}}
(\bs{M}, \mbox{---})\, $. 
Next, it is convenient to have some notation 
for {\it graded} Ext-groups. Thus, we set
$$
\underline{\rm Ext}^{n}_{\bs{A}}
(\bs{M}, 
{\boldsymbol N}) := 
\bigoplus_{d \in  {\mathbb Z}} 
\mbox{Ext}^{n}_{\bs{A}}
(\bs{M}, {\boldsymbol N}(d))\ ; 
$$
then $\,\underline{\rm Ext}^{n}_{\bs{A}}
(\bs{M}, \mbox{---})$ for 
$ n \geq 1 $ are the right derived functors 
of $ \underline{\rm Ext}^{0}_{
\bs{A}}
(\bs{M}, \mbox{---}) := 
\underline{\rm Hom}_{\bs{A}}
(\bs{M}, \mbox{---})\, $. 
To clarify this definition observe 
(see \cite{NV}, Corollary~I.2.12) that 
\begin{equation}
\la{ungr}
\underline{\rm Ext}^{n}_{\bs{A}}
(\bs{M}, {\boldsymbol N}) = 
\mbox{Ext}^{n}_{\mbox{\scriptsize Mod}
(\bs{A})}(\bs{M}, \bs{N}) \quad 
\mbox{for all}\ n \geq 0 \ ,
\end{equation}
at least when $\, \bs{M} \,$
is finitely generated. (On the right hand side
of (\ref{ungr}) $ \bs{M} $ and $ \bs{N} $ 
are regarded as objects 
in the category $ \mbox{\sf Mod}(\bs{A}) $ 
of {\it ungraded} $ \bs{A}$-modules.)

Finally, we mention the following
natural isomorphisms
(of graded vector spaces):
$$
\underline{\rm Ext}^{n}_{\bs{A}}
(\bs{M}, {\boldsymbol N}(d)) 
\cong
\underline{\rm Ext}^{n}_{\bs{A}}
(\bs{M}(-d), {\boldsymbol N}) 
\cong
\underline{\rm Ext}^{n}_{\bs{A}}
(\bs{M}, {\boldsymbol N})(d)
$$
valid for all $ d \in {\mathbb Z} $ 
and for all
$ \bs{M}, {\boldsymbol N} \in
\mbox{\sf GrMod}(\bs{A})\,$.

\subsection{Projective Schemes}
\label{S22}
Let $ \bs{A} $ be a Noetherian connected
graded $ k$-algebra, and let 
$ \bs{M} $ be a graded right module over 
$ \bs{A}\, $. We say that
$ m \in \bs{M} $ is a {\it $ \tau$orsion
element} of $ \bs{M} $ if 
$ m A_{n} = 0 $ for 
$ n \gg 0 \,$. The $\tau$orsion 
elements form a (graded) submodule in 
$ \bs{M} \,$; we denote it 
by $ \bs{\tau}\bs{M} \,$.
Equivalently,
$\, \bs{\tau}\bs{M} \,$ 
is the sum of all finite dimensional (over $ k $)
submodules of $ \bs{M}\, $. 
In particular, if $ \bs{M} $
is finitely generated, so is 
$ \bs{\tau}\bs{M} $ (since we
assume $ \bs{A} $ to be
Noetherian), and hence
$ \dim_{k} \boldsymbol{\tau}\bs{M} < 
\infty $ in that case.

A module $ \bs{M} $ is called a
{\it $\tau$orsion module} if $ \boldsymbol{\tau}\bs{M} = 
\bs{M}\,$, and {\it $\tau$orsion-free} if  
$ \boldsymbol{\tau}\bs{M} = 0 \,$.
The full subcategory of 
$ \mbox{\sf GrMod}(\bs{A}) $
consisting of all $\tau$orsion modules
will be denoted by 
$ \mbox{\sf Tors}(\bs{A})\,$.
Similarly, we shall write 
$ \mbox{\sf tors}(\bs{A})\,$
for the full subcategory of 
$ \mbox{\sf grmod}(\bs{A}) $
consisting of finitely generated 
$\tau$orsion modules. As we observed above,
the latter are precisely the graded 
modules which have finite dimension as 
vector spaces over $ k $.

Since both 
$ \mbox{\sf Tors}(\bs{A})\,$
and 
$ \mbox{\sf tors}(\bs{A}) $ 
are {\it dense} subcategories (that is, 
in any short exact sequence 
$ 0 \to \bs{M}' \to \bs{M} 
\to  \bs{M}'' \to 0\, $ the module
$\, \bs{M} $ is $\tau$orsion if and only if
$ \bs{M}' $ and $ \bs{M}'' $
are $\tau$orsion), we may introduce the 
quotient categories 
$$
\mbox{\sf Tails}(\bs{A}) :=
\mbox{\sf GrMod}(\bs{A})/
\mbox{\sf Tors}(\bs{A})\ , \quad
\mbox{\sf tails}(\bs{A}) :=
\mbox{\sf grmod}(\bs{A})/
\mbox{\sf tors}(\bs{A})\ .
$$
These are both abelian categories (see \cite{G}, pp. 367--369),
the second being a full subcategory of the first;
they are equipped with the exact projection functor
$ \pi: \mbox{\sf GrMod}(\bs{A}) \to 
\mbox{\sf Tails}(\bs{A}) $ which sends all the
$\tau$orsion objects in $ 
\mbox{\sf GrMod}(\bs{A}) $ to zero and
is universal (among additive functors) with respect 
to this property. 
Throughout the paper we shall denote quotient objects by 
script letters; for example, if 
$\, \bs{M} \in \mbox{\sf GrMod}(\bs{A})\,$,
we write $\, {\mathcal M}  := \pi \bs{M} \,$
for the corresponding object in 
$\,\mbox{\sf Tails}(\bs{A})\,$.
The shift in grading $ \bs{M} \mapsto \bs{M}(1) $
preserves $\tau$orsion modules,
hence carries over as an operation on 
quotient objects. The induced functor
$ {\mathcal M} \mapsto {\mathcal M}(1) $ on
$ \mbox{\sf Tails}(\bs{A}) $ (or on
$ \mbox{\sf tails}(\bs{A}) $) is called the
{\it twist functor}. 

In general, the description of the morphisms in 
$\,\mbox{\sf Tails}(\bs{A})\,$ is somewhat complicated. 
However, if $ \bs{M} $ is finitely generated, we have simply
\begin{equation}
\la{mor}
\mbox{\rm Hom}_{\,\mbox{\scriptsize \rm Tails}(\bs{A})\,}
({\mathcal M}, {\mathcal N}) \cong
\lim_{\longrightarrow}\, \mbox{\rm Hom}_{\bs{A}}
(\bs{M}_{\geq n}, {\boldsymbol N})\ ,
\end{equation}
where the system
$\, \{ \mbox{\rm Hom}_{\bs{A}}
(\bs{M}_{\geq n}, {\boldsymbol N})\} \,$ is 
directed by restriction of graded 
homomorphisms. It is easy 
to deduce from (\ref{mor}) when two objects in 
$ \mbox{\sf tails}(\bs{A}) $ are isomorphic, 
namely
\begin{equation}
\la{mor1}
{\mathcal M} \cong {\mathcal N}\ \mbox{in}\
\mbox{\sf tails}(\bs{A})\quad \Longleftrightarrow
\quad \bs{M}_{\geq n} \cong 
{\boldsymbol N}_{\geq n} \ \mbox{in}\ 
\mbox{\sf grmod}(\bs{A}) \ 
\mbox{for some}\  n \ .
\end{equation}
This perhaps explains the use of the word ``tails''.

If the algebra $\, \bs{A}\,$ is commutative and generated
by elements of degree one, then Serre's result tells us
that the categories $\, \Ta\,$ and $\, \ta \,$ are equivalent to the
categories of quasicoherent and coherent sheaves on the projective
scheme $\, X = \mbox{\sf proj}(\bs{A})\,$. For psychological reasons,
it is very helpful to use similar language also in the noncommutative 
case, even though in that case we shall not attempt to give any independent
meaning to ``$\mbox{\sf proj}(\bs{A})$''.
In what follows we shall refer to the objects of  
$ \mbox{\sf tails}(\bs{A}) $ 
(respectively, $ \mbox{\sf Tails}(\bs{A}) $)
as coherent (respectively, quasicoherent)
sheaves on $\, X = \mbox{\sf proj}(\bs{A})\,$, 
even when $ \bs{A} $ is not commutative, and we shall use 
the notation $\, {\mathcal O}_{X} := \pi \bs{A} \,$,
$  \mbox{\sf coh}(X) := 
\mbox{\sf tails}(\bs{A}) $,
$ \mbox{\sf Qcoh}(X) := 
\mbox{\sf Tails}(\bs{A}) \,$.

\subsection{Cohomology}
\label{S23}

In this section we outline the cohomology 
theory of coherent sheaves over noncommutative 
schemes, confining ourselves to results that 
will be used in the main body of the paper.
We keep the assumption that $ \bs{A} $ 
is a graded connected Noetherian $ k$-algebra.

For each $ {\mathcal M} \in 
\mbox{\sf Tails}(\bs{A}) $ the functor
$ \mbox{Hom}_{\mbox{\scriptsize \rm Tails}}({\mathcal M}, \mbox{---}) $
is left exact; since $ \mbox{\sf Tails}(\bs{A}) $ has 
enough injectives (see \cite{AZ}), the right derived functors
$ \mbox{Ext}^{n}({\mathcal M},\mbox{---}) $ are well defined.
As in the case of graded modules, we introduce the notation
$$
\underline{\rm Ext}^{n}({\mathcal M},{\mathcal N}) :=
\bigoplus_{d \in {\mathbb Z}} 
\mbox{Ext}^{n}({\mathcal M},{\mathcal N}(d))\ .
$$
\begin{definition}
\la{defcoh}
Let $\, {\mathcal M} \in \mbox{\sf Qcoh}(X) \equiv
\mbox{\sf Tails}(\bs{A}) \,$ be a 
quasicoherent sheaf over $\, X = \mbox{\sf proj}(\bs{A}) \,$. 
For each 
$\, n \geq 0 \,$  we define the {\it cohomology
groups} of $\, {\mathcal M} \,$ by 
\begin{equation}
\la{coh}
H^{n}(X, {\mathcal M}) := 
\mbox{Ext}^{n}({\mathcal O}_{X},{\mathcal M})\ ,
\end{equation}
where $ {\mathcal O}_{X} :=
\pi \bs{A} \in \mbox{\sf Qcoh}(X)\,$.
\end{definition}
The {\it cohomological dimension} of $ X $ is then defined by
$$
\mbox{\rm cdim}(X) := \mbox{\rm max}\{\, n \in {\mathbb N}\ :\ 
H^{n}(X,{\mathcal M}) \not= 0\ \mbox{for some}\  
{\mathcal M} \in \mbox{\sf Qcoh}(X) \,\}\ .
$$

Since $ \mbox{\sf Tails}(\bs{A}) $ is a 
$ k$-linear category, all the groups (\ref{coh}) are 
vector spaces over $ k $.  The graded objects
$$
\underline{H}^{n}(X, {\mathcal M}) := 
\bigoplus_{d \in {\mathbb Z}} 
H^{n}(X, {\mathcal M}(d))\ 
$$
are naturally graded right modules over 
$ \bs{A} \,$; we refer to them as the 
{\it full cohomology modules} of $ {\mathcal M} \,$. 
By (\ref{mor}), if $\, \ms{M} = \pi \bs{M}\,$, 
we have
$$
\underline{H}^{0}(X,\, \ms{M}) \cong 
\lim_{\longrightarrow} \, \underline{\rm Hom}_{\bs{A}}
(\bs{A}_{\geq n}, \bs{M})\ .
$$
From this it is easy to see that the functor 
$\, \omega :=  \underline{H}^{0}(X,\, \mbox{---}\,) \,$
is right adjoint to the projection functor $\, \pi:\, 
\Gma \to \Ta \,$. For any object $\, \ms{M} \in \Ta\,$,
the adjunction map $\, \ms{M} \to \pi\omega\ms{M}  \,$
is an isomorphism; the other adjunction map
$\, \bs{M} \to \omega \pi\bs{M} \cong \underline{H}^{0}(X, \ms{M}) \,$
is the restriction
\begin{equation}
\la{op}
\bs{M} \cong 
\underline{\rm Hom}_{\bs{A}}(\bs{A},\bs{M}) \to 
\lim_{\longrightarrow}\,\underline{\rm Hom}_{\bs{A}}
(\bs{A}_{\geq n}, \bs{M}) \cong \underline{H}^{0}(X, \ms{M})\ .
\end{equation}
Clearly, the kernel
of this map is $\, \bs{\tau M} \,$; more generally
(see \cite{AZ}, Prop.~7.2), there is an exact sequence 
of graded modules
\begin{equation}
\la{qqq}
0 \to \bs{\tau}\bs{M} \to \bs{M} \to \underline{H}^{0}(X, \ms{M}) \to
\lim_{\longrightarrow}\,\underline{\rm Ext}^{1}_{\bs{A}}
(\bs{A}/\bs{A}_{\geq n}, \bs{M}) \to 0\ .
\end{equation}

We shall need noncommutative versions of the 
basic theorems of Serre (finiteness and vanishing of 
cohomology, and duality). For this we have to impose
some additional conditions on our algebra $ \bs{A}\,$.
The least restrictive condition used in the literature
is the so-called $ \chi $-{\it condition} (see \cite{AZ}):
$\, \dim_{k} \underline{\rm Ext}^{n}_{\bs{A}}
({\boldsymbol k}_{\bs{A}}, \bs{M}) < \infty \,$ 
for all $\, \bs{M} \in \gma \,$ and for all $\, n \geq 0 \,$.
(Here $\, \bs{k}_{\bs{A}} := \bs{A}/\bs{A}_{\geq 1} \,$ 
denotes the ``trivial'' (right) module over $ \bs{A}\,$.)
The algebras occurring in the present paper
satisfy a much stronger condition: they are {\it Artin-Schelter} 
algebras. We shall concentrate on that case. 
The definition is as follows (see \cite{AS}).
\begin{definition}
\la{ASreg}
A graded connected algebra $ \bs{A} $
is called {\it Artin-Schelter} (or {\it Artin-Schelter regular}) 
if $ \bs{A} $ has 

$\,(i)$\ finite global dimension, say  
$\, \mbox{gl.dim}(\bs{A}) = d \,$;

$\,(ii)$\ polynomial growth, that is,
$\, \dim_{k}A_m \leq \gamma m^p $ for some
positive $ p \in {\mathbb Z} $, $ \gamma \in {\mathbb R} $, 
and  for all $ m \geq 0 \,$;

$\, (iii)$\  the (graded) {\it Gorenstein property}\,:
$\ \underline{\rm Ext}^{i}_{\bs{A}}
(\bs{k}_{\bs{A}}, \bs{A}) = 0 \,$ for all 
$ i \not= d $  and  $ \underline{\rm Ext}^{d}_{\bs{A}}
( \bs{k}_{\bs{A}}, \bs{A}) \cong 
\bs{k}_{\bs{A}}(l)\,$ for some integer $ l\,$ (called
the {\it Gorenstein parameter} of $\,\bs{A}\,$).
\end{definition}

If $ \bs{A} $ is commutative, then the
condition $ (i) $ in definition~\ref{ASreg}
already implies that $\, \bs{A} \,$ is isomorphic 
to a polynomial ring $\, k[x_0, x_1, \ldots, x_n] \,$
with some positive grading (see \cite{SZ}). 
Thus the only commutative Artin-Schelter 
algebras are polynomial algebras.
However, in the noncommutative case
there are many interesting examples (see \cite{AS}, 
\cite{ATV}, \cite{Steph} and references therein). 
The projective schemes associated with regular 
noncommutative algebras are referred to 
as ``quantum projective spaces''. The next theorem
provides further justification for this terminology.
\begin{theorem}[\cite{AZ}, Theorem~8.1]
\la{ASL}
Let $ \bs{A} $ be a Noetherian
Artin-Schelter algebra of global dimension 
$ d = n+1 $, and let $ X = 
\mbox{\sf proj}(\bs{A}) \,$.
Then $\,\mbox{\rm cdim}(X) = n \,$, and 
the full cohomology modules of 
$ {\mathcal O}_{X} := \pi \bs{A} $ are 
given by
$$
\underline{H}^{i}(X, 
{\mathcal O}_{X}) \cong
\left\{
\begin{array}{ccc}
   \bs{A}     & \quad \mbox{if} &  i = 0 \\
         0             & \quad \mbox{if} &  i \not= 0,\, n\\
\bs{A}^{*}(l) & \quad \mbox{if} & i = n
\end{array}
\right. \ ,
$$
where $ l $ is the Gorenstein parameter of $ \bs{A}\,$, and
$\, \bs{A}^{*} \,$ denotes the 
{\it graded dual} of $\, \bs{A} \,$ with components
$\, A_{i}^{*} := \mbox{\rm Hom}_{k}(A_{-i}, k)\,$.
\end{theorem}

Now we are in position to state the version of
Serre's theorems that we shall use.
\begin{theorem}
\la{FVT}
Let $ \bs{A} $ be a Noetherian Artin-Schelter 
algebra of global dimension $\, d = n+1 \,$ and Gorenstein 
parameter $ l\,$. Let $ X = 
\mbox{\sf proj}(\bs{A}) \,$.
Then if $\, {\mathcal M} \in \mbox{\sf coh}(X) \,$,
we have

$\,(a)$ \ \mbox{\rm (Finiteness)}  
$\ \dim_{k} H^{i}(X,{\mathcal M}) < \infty \,$ for all 
$\, i \geq 0 \,$; 

$\,(b)$ \ \mbox{\rm (Vanishing)}\  if $\, i \geq 1 \,$, then   
$\, H^{i}(X,\ms{M}(k)) = 0 \,$ for all $\, k \gg 0\,$;

$\,(c)$\ \mbox{\rm (Grothendieck-Serre Duality)}\ there are
natural isomorphisms
$$
\mbox{\rm Ext}^{i}({\mathcal M},
\ms{O}_{X}(-l)) \cong H^{n-i}(X, \ms{M})^{*}\ ,
$$
for $\, i = 0,1,2, \ldots, n\,$.
\end{theorem}
The Finiteness and Vanishing theorems have been proved by Artin and 
Zhang (see \cite{AZ}, Theorem~7.4) for any connected 
Noetherian algebra satisfying the $ \chi$-condition. 
For the Duality theorem we refer to \cite{YZ} (see also \cite{J}): 
we shall use this Theorem only in Section~\ref{Sect10}.

\subsection{\v{C}ech Cohomology}
\la{NCC}
In the classical case when $\, \bs{A} \,$ is commutative and generated by
elements of degree $ 1 \,$, we can calculate sheaf cohomology of the 
projective scheme $\, X = \pr{\bs{A}}\,$ using the \v{C}ech complex of
any affine open covering of $\, X \,$. Restricting a 
(quasicoherent) sheaf to an affine open set corresponds under Serre's
equivalence to (graded) localization of the associated 
$ \bs{A}$-module.
A natural noncommutative generalization of this construction 
has been suggested recently in \cite{Ver}, \cite{VW1}, \cite{VW2}.
Since in the noncommutative case one can define localization
only with respect to Ore sets, the existence of ``sufficiently many'' 
(homogeneous) Ore sets is a necessary condition to be imposed on the 
corresponding graded algebras. The class of such algebras 
(called {\it schematic}) is fairly rich and contains many 
interesting examples (see \cite{VW3}). Translating the definition
of a covering into algebraic language, we arrive at the following
\begin{definition}[see \cite{VW1}]
\la{scheme}
A Noetherian graded connected $k$-algebra $ \bs{A} $
is called {\it schematic} if there exists a 
finite number of (two-sided) homogeneous
Ore sets $\, U_1, U_2, \ldots,  U_s \,$ in $\, \bs{A} $ such that

$ (i)\ $ each $ U_{i} $ contains $ 1 \,$, and all
$\, U_{i} \cap \bs{A}_{\geq 1} \, $ are non-empty;

$ (ii)\ $ for any collection of elements $\, (u_1, u_2, \ldots, u_s) \in 
U_1 \times U_2 \times \ldots \times U_s \,$, there is an 
$\, m \in \mathbb{N} \,$ such that
\begin{equation}
\la{1}
\bs{A}_{\geq m} \subseteq \sum\limits_{i=1}^{s} u_i \bs{A} \ .
\end{equation}
A collection of Ore sets satisfying the conditions 
above is called a {\it covering of} $\, \bs{A} \,$.
\end{definition}
If $\, \bs{A} \,$ is schematic, let $\, N \,$ denote 
the least possible number of Ore sets 
covering $\, \bs{A} \,$. Following \cite{W}, we 
define the {\it schematic dimension} of $ \bs{A} $ by
$\,\mbox{sdim}(\bs{A}) := N-1 \,$.
If $ \bs{A} $ is commutative, $ \mbox{sdim}(\bs{A}) $ coincides 
with the usual (Krull or cohomological) dimension of the scheme 
$\, X = \mbox{proj}(\bs{A})\,$.
However, in the noncommutative case the schematic dimension
may happen to be smaller than $\,\mbox{cdim}(X)\,$, even for
Artin-Schelter algebras (see \cite{W} and Lemma~\ref{L4.1} below). 

The \v{C}ech complex of a covering of $\, \bs{A}\,$ 
is constructed in more or 
less the usual way, except that an ``intersection of open sets''
now depends on the order in which the sets intersect.
Fix a (finite) covering of $\, \bs{A} \,$, say 
$ \mathfrak{U} = \{  U_1, U_2, \ldots,  U_s \}\,$.
Given $ \bs{M} \in \Gma $ and a $ (p+1)$-tuple $ (i_{0},
i_{1}, \ldots, i_{p}) $ of indices, each $ i_{k} $ being
in $ \{1,2, \ldots, s \} $, we write
\begin{equation}
\la{8} 
\bs{M}_{i_{0} i_{1} \ldots i_{p}} := 
\bs{M} \underline{\otimes}_{\bs{A}} \bs{A}_{U_{i_{0}}} 
\underline{\otimes}_{\bs{A}} \ldots  \underline{\otimes}_{\bs{A}} 
\bs{A}_{U_{i_{p}}}\ ,
\end{equation}
where $ \bs{A}_{U_{i_{k}}} := \bs{A}[U_{i_{k}}^{-1}] $ is the 
(graded) localization of $ \bs{A} $ at $ U_{i_{k}} $. Now, 
for each $ p = 0, 1, 2, \ldots \,$, set
$$
\bs{C}^{p}(\mathfrak{U}, \ms{M}) :=
\bigoplus\limits_{(i_{0},
i_{1}, \ldots, i_{p})}^{} 
\bs{M}_{i_{0} i_{1}\ldots i_{p}} \, \in\, \Gma\ .
$$
Then $\, \bs{C}^{p}(\mathfrak{U}, \ms{M}) \,$ form a complex
of graded $ \bs{A}$-modules
\begin{equation}
\la{10} 
\bs{C}^{\bullet}(\mathfrak{U}, \ms{M}):\ 
0 \to \bs{C}^{0}(\mathfrak{U}, \ms{M})
\stackrel{d^{\,0}}{\longrightarrow} 
\bs{C}^{1}(\mathfrak{U}, \ms{M}) \stackrel{d^{\,1}}{\longrightarrow} \ldots 
\end{equation}
with coboundary maps $ d^{\,p}: \bs{C}^{p}(\mathfrak{U}, \ms{M}) \to 
\bs{C}^{p+1}(\mathfrak{U}, \ms{M}) $ defined in the usual way. 
For example, 
$$ 
d^{\,0}: \bigoplus_{i=1}^{s} \bs{M}_{i} \to \bigoplus_{i,j=1}^{s} 
\bs{M}_{ij}
$$
is given by the formula
$$
 d^{\,0}(m_1 \otimes u^{-1}_{1}, \ldots, m_s \otimes u_{s}^{-1})_{ij}
= m_{i} \otimes u_{i}^{-1} \otimes 1 - m_{j} \otimes 1 \otimes 
u_{j}^{-1} 
$$
in $ \bs{M}_{ij} = 
\bs{M} \underline{\otimes}_{\bs{A}} \bs{A}_{U_{i}}
\underline{\otimes}_{\bs{A}} \bs{A}_{U_{j}} \,$.

\vspace*{0.1cm}

The cohomology of the complex (\ref{10})
$$
\underline{\check{H}}{}^{p}(\mathfrak{U}, \ms{M}) := 
h^{p}[\bs{C}^{\bullet}(\mathfrak{U}, \ms{M})]\ ,
\quad p = 0, 1, 2, \ldots
$$
is called the (full) {\it \v{C}ech cohomology} of the sheaf $ \ms{M} $
relative to the covering $ \mathfrak{U}\,$. As in the commutative case,
we have the following general result.
\begin{theorem}
\la{Cech}
For all quasicoherent sheaves $  \ms{M} = \pi\bs{M} \in  \Qc{X} \,$,
and for any covering $ \mathfrak{U} $ of $ \bs{A} $,
there are natural isomorphisms of graded modules
$$
\underline{H}^{p}(X,  \ms{M}) \cong
\underline{\check{H}}{}^{p}(\mathfrak{U}, \ms{M})\ , 
\quad p = 0, 1, 2, \ldots \,.
$$
\end{theorem}
For the proof we refer the reader to \cite{VW2}.
We shall use Theorem~\ref{Cech} only for $\, p = 0 \,$, in which case it is
an elementary exercise. The isomorphism 
$\, \underline{H}^{0}(X,  \ms{M}) \to 
\underline{\check{H}}{}^{0}(\mathfrak{U}, \ms{M}) \,$ is defined as 
follows. Recall that an element of 
$\,\underline{H}^{0}(X, \ms{M}) \,$ is represented by a 
homomorphism $\, f : \bs{A}_{\geq n} \to \bs{M} \,$ (for some $\, n $).
Condition $\,(i)\,$ in Definition~\ref{scheme} implies that
$\, (\bs{A}_{\geq n})_{U_i} \cong \bs{A}_{U_i} \,$ for any $\, n\,$;
so after localization $\, f \,$ defines a homomorphism 
$\, f_{i}: \bs{A}_{U_i} \to \bs{M}_{U_i} \,$ for each $ i\,$.
Assigning to $\, f \,$ the element 
$\, (f_{1}(1),\, f_{2}(1),\, \ldots, \,f_{s}(1)) 
\in \bs{C}^{0}(\mathfrak{U}, \ms{M})\,$
defines the desired isomorphism. Combining this with (\ref{op})
we get
\begin{proposition}
\la{ch}
Let $\, \bs{M} \,$ be a graded $ \bs{A}$-module, 
$\, \ms{M} = \pi\bs{M} \,$ the associated sheaf.
Then the natural map
$\,\bs{M} \to \underline{H}^{0}(X,\ms{M}) \to 
\underline{\check{H}}{}^{0}(\mathfrak{U},\ms{M}) \,$
is given by 
$$ 
m  \mapsto  (m \otimes 1, \ldots, m \otimes 1)\ \in
\  \bs{C}^{0}(\mathfrak{U}, \ms{M}) = 
\bs{M}_{U_{1}} \oplus \ldots \oplus \bs{M}_{U_{s}}\ .
$$
\end{proposition}

We shall denote the $ d$-th graded component of  
$\,\underline{\check{H}}{}^{p}(\mathfrak{U},\ms{M}) \,$
by $\, \check{H}^{p}(\mathfrak{U},\ms{M}(d)) \,$.
We then have
$$
\underline{\check{H}}{}^{p}(\mathfrak{U},\ms{M})
= \bigoplus_{d \in {\mathbb Z}} 
\check{H}^{p}(\mathfrak{U},\ms{M}(d)) \ ,
$$
in conformity with our usual notation.

\section{The Weyl Algebra and its Homogenizations} 
\la{Sect3}

From now on, we set $\, k = {\mathbb C} \, $. Let 
$\, A = {\mathbb C}\langle x,y \rangle/([x,y]-1) \,$ be
the first Weyl algebra over $ {\mathbb C}\, $. 
Unlike the commutative algebra $\, {\mathbb C}[x,y]\, $,
the Weyl algebra admits no natural grading; 
however, it has many natural filtrations.

\subsection{Filtered Rings and Modules} 
\label{SFr}
We recall that a {\it filtration} on an algebra
$ A $ is an increasing sequence of linear subspaces
\begin{equation}
\la{filt}
\ldots\, \subseteq A_{k-1} \subseteq A_{k} \subseteq A_{k+1}\, \ldots \ ,
\end{equation}
indexed by the integers, such that $\, 1 \in A_{0} \,$,
$\, \bigcup_{k \in {\mathbb Z}} A_{k} = A \, $ and 
$\, A_{j} A_{k} \subseteq A_{j+k} \,$ for all
$\, j, k \in {\mathbb Z} \,$.
A filtration is called {\it positive} if $ A_{k} = 0 $ for all 
$ k < 0 \,$. If $ A $ is a filtered algebra, 
then a (right) $ A$-module $ M $ is called a 
{\it filtered $ A$-module} if there is an increasing 
sequence of linear subspaces
\begin{equation}
\la{5.0}
\ldots \subseteq M_{k-1} \subseteq M_{k}
\subseteq M_{k+1} \subseteq \ldots \ ,
\end{equation}
such that $\, 
\bigcup_{k \in {\mathbb Z}} M_k  = M\, $
and $\, M_k A_j \subseteq M_{k+j} \,$ for
all $\, k, j \in {\mathbb Z}\,$. We shall
assume the filtration (\ref{5.0}) to be 
{\it separated}, meaning that 
$ \bigcap_{k \in {\mathbb Z}} M_k  = 0 \,$.

Attached naturally to a filtered algebra are the following two 
graded algebras:
$$
\bs{A} := 
\bigoplus\limits^{}_{k \in  {\mathbb Z}} A_{k}\quad  , \quad
\boldsymbol{GA} := 
\bigoplus\limits^{}_{k \in  {\mathbb Z}} A_{k}/A_{k-1}\ .
$$
The algebra $ \bs{A} $ is called the {\it Rees algebra} 
of $ A $ with respect to the filtration (\ref{filt}). 
It can be identified with a subring of the ring of 
Laurent polynomials (in one variable $ t $) 
with coefficients in $ A $. To be precise, we have
\begin{equation}
\la{Laur}
\bs{A} \cong 
\bigoplus\limits^{}_{k \in \mathbb{Z}} A_{k} t^k 
\hookrightarrow A[t, t^{-1}] \ ,
\end{equation}
where the grading on $ A[t, t^{-1}] $ is defined by
$ \deg(t) = 1 $ and $ \deg(a) = 0 $ for all elements 
$ a \in A \,$. The  algebra $ \bs{GA} $ is called
the {\it associated graded ring} of $ A \,$.
Under the identification (\ref{Laur}) there is
a natural isomorphism of graded algebras
\begin{equation}
\la{GA}
\boldsymbol{GA} \cong \bs{A}/\langle t\rangle\ ,
\end{equation}
where $ \langle t \rangle $ denotes the two-sided ideal
of $ \bs{A} $ generated by the central element $\, t \,$. 

Similarly, if $ M $ is a filtered $ A$-module
then we have the {\it Rees module} 
\begin{equation}
\label{RM}
\bs{M} :=
\bigoplus\limits^{}_{k \in  {\mathbb Z}} M_{k}\ 
\in\ \mbox{\sf GrMod}(\bs{A})\ ,
\end{equation}
and the {\it associated graded module} 
$$
\bs{GM} :=
\bigoplus\limits^{}_{k \in  {\mathbb Z}} M_{k}/M_{k-1}\ 
\in\ \mbox{\sf GrMod}(\boldsymbol{GA})\ .
$$
Identifying $ \bs{A} $ with a ring of 
$ A$-valued Laurent polynomials 
(see (\ref{Laur})), we have $\, \bs{M} \cong 
\bigoplus\limits^{}_{k \in  {\mathbb Z}} M_{k} t^k 
\hookrightarrow M[t, t^{-1}] \, $, where
$ \, M[t, t^{-1}] := M \otimes_{A} A[t, t^{-1}] \,$,
and hence, in view of (\ref{GA}), the following 
isomorphisms of graded $ \bs{GA}$-modules
\begin{equation}
\la{Laur1}
\boldsymbol{GM} \cong \bs{M}/\bs{M}t 
\cong \bs{M} \otimes_{\bs{A}} 
\bs{A}/ \langle t \rangle\ 
\cong \bs{M} \otimes_{\bs{A}} 
\bs{GA}\ .
\end{equation}

When $ A $ is commutative, the above constructions
have a simple geometrical meaning: 
$\, X := \pr{\bs{A}} \,$ is a projective scheme 
containing the affine scheme $\, \mbox{\sf Spec}(A) \,$
as an open subset, and the sheaf 
$\, \ms{M} := \pi \bs{M} \,$ is an extension to $ X $
of the sheaf $ \tilde M $ on $\, \mbox{\sf Spec}(A) \,$ 
corresponding to $\, M \,$. Thus, from an algebraic 
point of view, the projective compactification $ X \,$
is determined by the choice of filtration on $ A \,$, 
and the extension of $ \tilde M $ to $ X\,$ is then determined 
by the choice of filtration on $ M \,$. Furthermore,
$\, \pr{\bs{GA}} \,$ is the ``hypersurface at infinity''
(in our case it will be a line) in $ X\,$, and 
$\, \pi \bs{GM} \,$ is the restriction of $\, \ms{M}\,$
to this hypersurface. We shall use similar language also
in the noncommutative case.
\subsection{Weight Filtrations} 
\label{SBf}

We now introduce the class of filtrations on the Weyl algebra
$ A $ which we shall use in the present paper. 
Given a pair of positive integers $\, {\boldsymbol w} := 
(w_1, w_2)\, $, we set
\begin{equation}
\la{Bf}
A_{k}({\boldsymbol w}) := 
\mbox{span}_{{\mathbb C}}
\{\, x^{\alpha} y^{\beta}\, | \, w_1 \alpha + w_2 
\beta \leq k\,
\} \subset A
\end{equation}
for each $ k \in  {\mathbb Z}\, $. Then $\, 
\{ A_{\bullet}({\boldsymbol w}) \} $ is a positive 
locally finite filtration on 
$ A $ with $ A_{0}({\boldsymbol w}) = {\mathbb C} \,$.
We call (\ref{Bf}) the {\it filtration of weight} 
$ {\boldsymbol w} \,$. In particular,
if $ {\boldsymbol w} = (1,1) $, this is the standard
Bernstein filtration on $ A \,$. 

To describe the Rees algebra associated with 
(\ref{Bf}) we use the identification (\ref{Laur}).
Setting $ X := x  \cdot t^{w_1} $, $ Y := y  \cdot t^{w_2} $ 
and $ Z := 1 \cdot t \,$, we observe that $ \bs{A} $ 
is isomorphic to the graded algebra generated 
(over $ {\mathbb C} $) by 3 elements $ X, Y $ and
$ Z $ (in degrees $ w_1 $, $ w_2 $ and $ 1 $ respectively)
subject to the defining relations
\begin{eqnarray}
\la{4.4}
XZ - ZX & = & 0 \ , \nonumber \\
YZ - ZY & = & 0\ , \\
XY - YX & = & Z^{|{\boldsymbol w}|} \ , \nonumber
\end{eqnarray}
where $ |{\boldsymbol w}| := w_1 + w_2 \,$.
We call this algebra
the {\it homogenized Weyl algebra of weight} 
$\, {\boldsymbol w} \,$ and denote it
by $ \bs{A}({\boldsymbol w})\,$ 
(or simply by $\,\bs{A} \,$ when there is 
no danger of confusion).

The following two propositions collect some basic 
properties of $ \bs{A}({\boldsymbol w}) $.
\begin{proposition}
\la{L4.1}
For every positive weight vector, the algebras 
$\, \bs{GA}({\boldsymbol w}) \,$ and 
$\, \bs{A}({\boldsymbol w}) \,$ 
are Noetherian Artin-Schelter algebras of global 
dimensions 2 and 3 respectively. The corresponding
Gorenstein parameters are $ |\bs{w}| $ and $ |\bs{w}| + 1\,$.
\end{proposition}

\begin{proof} We have
\begin{equation}
\la{4.7}
\bs{GA}({\boldsymbol w}) \cong 
\bs{A}({\boldsymbol w})/ \langle Z \rangle  \cong 
{\boldsymbol S}({\boldsymbol w})\ ,
\end{equation}
where $\, \bs{S}(\bs{w}) := {\mathbb C}[\bar{x}, \bar{y}] \,$
is the graded {\it commutative} polynomial ring in two variables of 
weight $ \bs{w} \,$ (the first isomorphism in (\ref{4.7}) is just 
(\ref{GA}), while the second follows immediately from the defining 
relations (\ref{4.4}).) Hence $\,\bs{GA}(\bs{w})\,$ has the properties
stated in the lemma. Since $\, \bs{GA}(\bs{w}) \,$
is Noetherian and has global dimension 2, 
$ \bs{A}(\bs{w}) $ is also Noetherian and has global dimension 
3 (see \cite{L}, Proposition~3.5 and \cite{LO}, Theorem~II.8.2
respectively). That $ \bs{A} $ is Artin-Schelter
follows from \cite{L}, Theorem~5.10 and Theorem~6.3. 
According to \cite{ATV}, Proposition~2.14, the Gorenstein 
parameter $ l $ is equal to the degree of the inverse
of the Poincar\'e series 
$\, P_{\bs{A}}(s) := \sum_{k \geq 0}^{}\, 
\dim_{\mathbb C} A_{k}(\bs{w})\,s^k\,$. Since $ \bs{A} $ is isomorphic 
(as a graded vector space)
to the commutative polynomial ring in three variables 
of weights $ (w_1, w_2, 1)\,$, we have
$$
P_{\bs{A}}(s) = \frac{1}{(1-s^{w_1})(1-s^{w_2})(1-s)}\ ,
$$
and therefore the Gorenstein parameter is $\, w_1 + w_2 +1\,$.
\end{proof}
\begin{proposition}
\la{L4.11}
For every positive weight vector, 
$\, \bs{A}({\boldsymbol w}) \,$ is a schematic algebra of 
schematic dimension $ 1\,$.
\end{proposition}
\begin{proof} 
This is proved in \cite{W} in the case $\,\bs{w} = (1,1)\,$; 
more precisely, it is shown in \cite{W} that $ \bs{A} $ can be covered by
the two Ore sets consisting of the powers of $ X $ and of $ Y $. 
The proof in general is similar, but we shall work with a ``finer''
covering, that is, with  larger Ore sets. This covering will be needed 
in Section~\ref{Sect6}.

For an element $\, q(x,y) \in A \,$ we denote by 
$\, \bs{q}(X,Y,Z) \in \bs{A} \,$ 
its homogenization in $ \bs{A} $, that is, we set $ \bs{q}(X,Y,Z) = 
Z^{d} q(X/Z^{w_1}, Y/Z^{w_2}) \,$, where $ d = \deg_{\bs{w}}(q)\,$. 
Now define
\begin{equation}
\la{4}
\begin{array}{rcl}
U_{1} &:=& \{\, \bs{q}_1(X,Z) \in \bs{A}\,|\, q_{1}(x) \in 
\c[x] \setminus \{0\}\,\}\ ,\\*[1ex]
U_{2} &:=& 
\{\, \bs{q}_2(Y,Z) \in \bs{A}\, |\, q_{2}(y) \in 
\c[y] \setminus \{0\}\,\} \ .
\end{array}
\end{equation}
We claim that $\,\mathfrak{U} := \{U_1,  U_2\} \,$ is a  
covering of $ \bs{A}\,$.  
Clearly, $\, U_1 \,$ and $\, U_2 \,$ are multiplicatively closed 
subsets of $\, \bs{A}\,$; and they satisfy the Ore condition (on both sides) 
because they consist of (homogeneous) 
locally ad-nilpotent elements in $ \bs{A}\,$.
On the other hand, for any (nonzero) polynomials
$ q_{1}(x) $ and $ q_{2}(y) $, the ideal $\, q_{1}(x)A +  q_{2}(y)A \,$ has 
finite codimension in $ A $
and hence coincides with $ A\,$.
It follows that $ 1 = q_{1}(x) a + q_{2}(y) b $
for some $ a,b \in A \,$, and therefore
$ Z^{r} \in \bs{q}_1 \bs{A} + \bs{q}_2 \bs{A} $
for some $ r > 0 $. If either $ q_{1} \in \mathbb{C} $  or
$ q_{2} \in \mathbb{C} \,$, the condition (\ref{1}) holds 
trivially. So we may assume
that $ q_{1} $ and $ q_{2} $ both have positive degree, say 
$ \deg_{\bs{w}}(q_{1}) = w_1 d_1 > 0 $ and 
$ \deg_{\bs{w}}(q_{2}) = w_2 d_2 > 0 \,$. Then 
$ \bs{q}_1(X,Z) - X^{d_1} \in 
Z \bs{A} $ and $ \bs{q}_2(Y,Z) - Y^{d_2} \in 
Z \bs{A} \,$, hence
$\, X^{k d_1} \in \bs{q}_1 \bs{A} + Z^k \bs{A}\, $  and  
$\, Y^{k d_2} \in \bs{q}_2 \bs{A} + Z^k \bs{A}\, $
for any $\, k \geq 1 \,$.  Taking $\, k = r\,$, we find that there is an 
$m$ such that
$$
\bs{A}_{\geq m} \subseteq  X^{r d_1} \bs{A}  + Y^{r d_2} \bs{A} + Z^r \bs{A} 
\subseteq  \bs{q}_1 \bs{A} + \bs{q}_2 \bs{A} \ ,
$$
which shows that the pair of Ore sets 
$ \{U_1,  U_2\} $ covers the 
algebra $ \bs{A} \,$, as required.
\end{proof}
\subsection*{Remark}\ 
Unlike global or Gel'fand-Kirillov dimension,
the schematic dimension distinguishes $ \bs{A} $
from the commutative polynomial algebra $\, \mathbb{C}[X,Y,Z] \,$.

\subsection{Weighted Projective Planes}
Given a weight vector $ \bs{w} \,$, we write 
$$ 
\P(\bs{w}) := \mbox{\sf proj}\, \bs{A}(\bs{w})\quad , \quad
\mathbb{P}^{1}(\bs{w}) := \mbox{\sf proj}\, \bs{GA}(\bs{w})
$$
for the (hypothetical) projective schemes associated to $ \bs{A} $ and $ \bs{GA} \,$.
The identification $\, \boldsymbol{GA} \cong 
\bs{A}/ \langle Z \rangle \, $ provides a natural
epimorphism of graded algebras
\begin{equation}
\la{4.10}
\boldsymbol{i} : \bs{A} \to \boldsymbol{GA} \ .
\end{equation}
As usual, we have the functors 
$\, \bs{i}_{*} \,$
and $\, \bs{i}^{*} \,$ of restriction and extension 
of scalars (if $ \bs{M} $ is a (right graded) 
$ \bs{GA}$-module, then $\, \bs{i}_{*} \bs{M} \,$ is the 
same vector space $ \bs{M} $ with $ \bs{A}$-module structure 
defined via (\ref{4.10}), while if $ \bs{M} $  is an 
$ \bs{A}$-module, $\, \bs{i}^{*} \bs{M} \,$ is the 
$ \bs{GA}$-module $\,\bs{M} \underline{\otimes}_{\bs{A}} \bs{GA}\,$).
These functors both preserve the classes of finitely 
generated and of $\tau$orsion 
modules, and hence descend to functors $ i_* $ and $ i^* $ 
on the categories of coherent (or quasicoherent) sheaves over
$ \P $ and $ \mathbb{P}^{1}\, $. We shall call $ \mathbb{P}^{1} $ the {\it line at infinity}
in $ \P $ and sometimes denote it by $ l_{\infty}\,$.
If $ \ms{M} $ is a coherent sheaf over $ \P \,$, we call
$ i^* \ms{M} $ the {\it restriction of $ \ms{M} $ to the line
at infinity}.

For future use we record the simple
\begin{lemma}
\la{seq}
For any $\, \ms{M} \in \C\P \,$ there is an exact sequence
$$
\ms{M}(-1) \to \ms{M} \to i_* i^*\ms{M} \to 0\ ,
$$
where the first map is induced by multiplication by
$\, Z \in \bs{A}(\bs{w})\,$.
\end{lemma}
\begin{proof}
If $ \bs{M} $ is a graded $ \bs{A}$-module with 
$\, \ms{M} = \pi\bs{M}\,$, it follows at once 
from (\ref{Laur1}) that the graded 
quotient $\, \bs{M}/\bs{M}Z \,$ is canonically 
isomorphic to $\, \bs{i}_* \bs{i}^*\bs{M}\,$.
Applying the (exact) functor $ \pi $ to 
$$ 
\bs{M}(-1) \stackrel{\cdot \, Z}{\longrightarrow} \bs{M} \to 
\bs{i}_* \bs{i}^*\bs{M} \to 0 
$$
we get the lemma.
\end{proof}
\subsection*{Remark}
Nearly all the results in this paper remain true if we replace
$\, \bs{A}(\bs{w}) \,$ by the {\it commutative} graded algebra 
$\, \bs{A}_{0} := \c[X,Y,Z] \,$ with weights $\,(w_1, w_2, 1)\,$.
However, except when $\, w_1 = w_2 = 1 \,$, we do not have
a Serre equivalence between the category of coherent sheaves 
(in the usual sense) over $\,\pr{\bs{A}_{0}}\,$ and 
the category $\, \mbox{\sf tails}(\bs{A}_{0}) \,$. 
Our results would refer to the latter category,
and so (probably) would not give much information about 
the usual ``weighted projective spaces'' studied in (for example)
\cite{Del}, \cite{Dol} and \cite{BR}. For a similar reason, 
our $\, \P(\bs{w}) \,$ are different from the quantum 
weighted projective planes introduced recently in \cite{St1}.

\section{The Linear Data Associated to an Ideal}
\la{Sect4}
Let $ M $ be a finitely generated torsion-free rank one right
module over the Weyl algebra $ A\,$. We fix a (positive)
filtration of weight $ \bs{w} $ on $ A \,$. We also fix,
temporarily, an embedding 
of  $ M $ as an ideal in $ A \,$; then we have the {\it induced 
filtration} $\, M_k = M \cap A_k \,$ on $ M\,$. It is easy to see
that, up to an overall shift, the filtration is independent of 
the choice of embedding: our first task is to normalize this 
overall shift.  Since $\, M \subseteq A \,$, the 
corresponding Rees module $\, \bs{M} \,$ (see (\ref{RM}))
is a graded ideal in $ \bs{A} \equiv \bs{A}(\bs{w})\,$.
Let $\, \ms{M} = \pi \bs{M} \,$  be the associated sheaf
over $ \P(\bs{w})\,$.
\begin{lemma}
\la{LL1}
There is a unique $\, a \in {\mathbb Z} \,$ such that
the restriction of $ \ms{M}(a) $ to the line at infinity
in $ \P $ is trivial. 
\end{lemma}
\begin{proof}
We have $\, \bs{i}^{*}\bs{M} = \bs{GM}\,$ (see (\ref{Laur1})).
The embedding of $\, \bs{M}\,$ in $\, \bs{A} \,$ induces an 
embedding of $\, \bs{GM} \,$ in  $\, \bs{GA} \,$ as a homogeneous
ideal. Now,  $\, \bs{GA} \,$ is just a commutative polynomial 
algebra in two variables; hence if $ f $ is the greatest common 
divisor of the elements of $\, \bs{GM} \,$, then  
$\, f^{-1} \bs{GM}  \,$ is a (homogeneous) ideal of 
finite codimension in $\,\bs{GA}\,$.
Denoting by $ a $ the degree of $ f $ in $ \bs{GA} \,$, we therefore
have an exact sequence of graded $ \boldsymbol{GA}$-modules 
\begin{equation}
\label{5.8}
0 \to \boldsymbol{GM}(a) \to \boldsymbol{GA} \to
\boldsymbol{GA}/\boldsymbol{GM}(a) \to 0 \,
\end{equation}
with finite-dimensional quotient term.  The quotient functor 
$\pi$ annihilates finite-dimensional modules, so applying 
$\, \pi \,$ to (\ref{5.8}), we get the desired 
isomorphism
$$
i^* \ms{M}(a) \cong \ms{O}_{\mathbb{P}^{1}} \quad \mbox{in}\ 
\C \mathbb{P}^{1}\ . 
$$
The uniqueness of $\, a \,$ follows from 
the fact that $ \ms{O}_{{\mathbb P}^{1}}(k) \cong  
{\mathcal O}_{{\mathbb P}^{1}} $ in 
$\,\mbox{\sf coh}({\mathbb P}^{1})\,$ 
only if $ k = 0\,$. Indeed, assuming the contrary, by (\ref{mor1})
we have
$\, \boldsymbol{GA}(k)_{\geq N} \cong 
\boldsymbol{GA}_{\geq N} \,$ for some $\,N\,$,
and therefore
$\, \dim_{\mathbb{C}} \boldsymbol{GA}_{n+k} = \dim_{\mathbb{C}}
\boldsymbol{GA}_{n} \,$ for all $ n \geq N \,$. This implies that 
the sequence of 
numbers $\,\left(\dim_{\mathbb{C}} \boldsymbol{GA}_{n}\right)\,$ 
is bounded,  which is obviously not the case.
\end{proof}
\begin{lemma}
\la{min}
Let $\, \delta \,$ be the minimum filtration degree of elements
of $ M \,$. Then $\,  \delta \geq a \,$; if $ M $ is not cyclic
then $\, \delta > a \,$.
\end{lemma}
\begin{proof}
As in the proof of the preceding Lemma, we identify $\, \bs{GM} \,$
with an ideal in the polynomial ring $\, \bs{GA} \,$; then
$\, \delta \,$ is the minimum degree of elements in $\, \bs{GM} \,$,
and hence $\, \delta \geq a\,$. If $ \delta = a \,$, then 
$\, \bs{GM} \,$ is cyclic (generated by the greatest common 
divisor $ f$ above), and hence $ M $ is also cyclic.
\end{proof}
\begin{proposition}
\la{hzero}
The natural map 
$\, \bs{M} \to \underline{H}^{0}(\P\,,\, \ms{M}) \,$
in (\ref{qqq}) is bijective.
\end{proposition}
\begin{proof}
It is obvious that $\, \bs{\tau M} = 0\,$, so we have only 
to prove that the Ext term in (\ref{qqq}) is zero. Let 
$\, \bs{N} := \bs{A}/\bs{M} \,$. We show first that
\begin{equation}
\la{tauN}
\lim_{\longrightarrow}\,\underline{\rm Ext}^{1}_{\bs{A}}
(\bs{A}/\bs{A}_{\geq n}, \bs{M}) \cong \bs{\tau N}.
\end{equation}
For brevity, set $\,\bs{A}_{<n} := \bs{A}/\bs{A}_{\geq n}\,$.
Applying the functor $\, \underline{\rm Hom}_{\bs{A}}
(\bs{A}_{<n},\,\mbox{---}) \,$ to 
$\,0 \to \bs{M} \to \bs{A} \to \bs{N} \to 0 \,$, we get
the exact sequence
$$
\underline{\mbox{Hom}}_{\bs{A}}
(\bs{A}_{<n}, \bs{A}) 
\to \underline{\mbox{Hom}}_{\bs{A}}
(\bs{A}_{<n}, {\boldsymbol N}) \to
\underline{\mbox{Ext}}^{1}_{\bs{A}}
(\bs{A}_{<n},\bs{M}) 
\to \underline{\mbox{Ext}}^{1}_{\bs{A}}
(\bs{A}_{<n}, \bs{A}) \ .
$$
The first term in this sequence is obviously zero; and the 
Gorenstein property (see Definition~\ref{ASreg}) of 
$ \bs{A} $ implies 
that the last term is also zero, because  $\, \bs{A}_{<n} \,$ has 
finite length. Thus $\, \underline{\mbox{Hom}}_{\bs{A}}
(\bs{A}_{<n}, {\boldsymbol N}) \cong 
\underline{\mbox{Ext}}^{1}_{\bs{A}}(\bs{A}_{<n},\bs{M})\,$
for all $ n \,$. Passing to the limit as $\, n \to \infty \,$,
we get (\ref{tauN}). Hence Proposition~\ref{hzero} follows if 
we show that  $\, \bs{\tau N} = 0 \,$. 
Suppose that $\, a \in A_k \,$ represents a $\tau$orsion element
in $\,\bs{N}\,$. This means that for some $ n \geq 0 $
we have $ a \bs{A}_{\geq n} \subset \bs{M} \,$, and hence
$\, a {A}_{n} \subset M_{n+k} \,$. Since $\, 1 \in A_{n}\,$,
we find that $\, a \in M_{n+k} \cap A_{k} = M_{k} \,$, and 
hence $ a $ represents zero in $\, \bs{N}\,$. Thus
$\, \bs{\tau N} = 0\,$.
\end{proof}

We use Lemma~\ref{LL1} to fix the ambiguous shift in the induced
filtration. If $\, M \,$ is an (embedded) ideal, then the 
uniqueness of the number $\, a \,$ in  Lemma~\ref{LL1} shows
that the filtration $\, M_{k}^{\circ} := M_{k+a}\,$ on $\, M \,$
is independent of the choice of embedding. We shall call
this filtration the {\it normalized} induced filtration 
on $\, M \,$. From Lemma~\ref{min}, we get 
\begin{lemma}
\la{min1}
Let $\, M\,$ be an ideal of $\, A \,$ with the normalized
induced filtration, and let $\, d \,$ be the minimum 
filtration degree of elements of $\, M\,$. Then 
$\, d \geq 0 \,$, and if $\, M \,$ is not cyclic,
then $\, d > 0 \,$.
\end{lemma}
From now on, changing notation, 
$\, \ms{M} \,$ will always denote the extension of $\, M\,$
to $\, \P(\bs{w})\,$ determined by the normalized
induced filtration (so that $\, \ms{M}|_{l_{\infty}}\,$
is trivial). We call $\, \ms{M} \,$ the 
{\it canonical extension} of $\, M\,$.  The next Theorem gathers 
together the information we need about the cohomology of $\,\ms{M}\,$.
\begin{theorem}
\la{PP1} 
Let $\,\ms{M}\,$ be the canonical extension of an ideal of $\,A\,$.  
Then

\smallskip
$ (i)\ $ 
The map $\, H^{1}(\P, \ms{M}(k-1)) \to  H^{1}(\P,\ms{M}(k)) \,$
induced by multiplication by $ Z $ is injective for 
$ k < 0 $ and surjective for $\, k > -|\bs{w}|\,$.

$ (ii)\ $ We have
\begin{eqnarray}
&& H^{0}(\P, \ms{M}(k)) = 0 \quad \mbox{for}\ k < 0\ ,\nonumber \\
&& H^{2}(\P, \ms{M}(k)) = 0 \quad \mbox{for}\ k \geq  -|\bs{w}|\ . 
\nonumber
\end{eqnarray}

$ (iii)\ $ Furthermore, if $ M $ is not cyclic, we have also
$\, H^{0}(\P, \ms{M}) = 0\,$, and 
$$
\dim_{\mathbb{C}} H^{1}(\P, \ms{M}) = \dim_{\mathbb{C}} 
H^{1}(\P, \ms{M}(-1)) - 1\ .
$$
\end{theorem}
\begin{proof}
The map of sheaves $\, \ms{M}(k-1) \to  \ms{M}(k)\,$
induced by multiplication by $ Z $ is clearly injective 
(for any $\, k \in \Z$). Using Lemma~\ref{seq} and bearing 
in mind that $\, i^* \ms{M} \cong \ms{O}_{\mathbb{P}^{1}} \,$ we get 
the short exact sequence
\begin{equation}
\la{5.3}
0 \to {\mathcal M}(k-1) \to {\mathcal M}(k) \to 
i_{*} \ms{O}_{\mathbb{P}^{1}}(k) \to 0\ .
\end{equation}
By Theorem~8.3 of \cite{AZ}, we have
$$
H^{i}(\P, \, i_{*} \ms{O}_{\mathbb{P}^{1}}(k)) \cong 
H^{i}(\mathbb{P}^{1},\, \ms{O}_{\mathbb{P}^{1}}(k))\quad 
\mbox{for all} \quad i \geq 0\ .
$$
On the other hand, by Theorem~\ref{ASL} and Proposition~\ref{L4.1},
we have
\begin{equation}
\la{5.14}
H^{i}(\mathbb{P}^{1},\, \ms{O}_{\mathbb{P}^{1}}(k)) = 
\left\{
\begin{array}{ccc}
S_{k}\,  & \quad \mbox{when} & i=0 \\*[1ex]
(S_{-k-|\bs{w}|})^*\,  & \quad 
\mbox{when} & i=1\\*[1ex]
0  & \quad \mbox{when} & i>1
\end{array}
\right. 
\end{equation}
where $ S_{k} $ is the $ k$-th graded component of the
commutative polynomial algebra $\, \bs{S}(\bs{w}) = 
\c[\bar x, \bar y] \,$ of weight $ \bs{w}\,$.
Therefore the first and last terms of the exact sequence
$$
H^{0}(\P, i_{*} \ms{O}_{\PP^{1}}(k)) \to  
H^{1}(\P, \ms{M}(k-1)) \to  H^{1}(\P,\ms{M}(k)) \to 
H^{1}(\P, i_{*} \ms{O}_{\PP^{1}}(k))
$$
coming from (\ref{5.3}) are isomorphic to $\, S_{k} \,$
and $\, (S_{-k-|\bs{w}|})^* \,$ respectively. Since the grading
on $ \bs{S}(\bs{w}) $ is positive, part $ (i) $ of the Theorem 
follows.

The assertion about $\, H^{0} \,$ in part $ (ii) $ is immediate 
in view of Proposition~\ref{hzero} and Lemma~\ref{min1}. 
To prove the assertion about
$\, H^{2} \,$ we observe (again looking at the long cohomology 
exact sequence of (\ref{5.3})) that the map
$$
H^{2}(\P, \ms{M}(k-1)) \to  H^{2}(\P,\ms{M}(k))
$$
is an isomorphism for $\, k-1 \geq -|\bs{w}|\,$. 
By the Vanishing Theorem~\ref{FVT}$(b)$, 
$\, H^{2}(\P,\ms{M}(k))\,$ is zero for $\, k\gg 0 \,$,
hence it is zero for all $\, k \geq -|\bs{w}|\,$.

It remains to prove part $ (iii) $ of the Theorem.
The fact that $\, H^{0}(\P, \ms{M}) = 0\,$ again follows 
 from Lemma~\ref{min1} and Proposition~\ref{hzero}.
From (\ref{5.3}) (with $\,k=0$) we get the exact sequence
$$
0 = H^{0}({\mathbb P}^{2}_{q}, {\mathcal M}) 
\to {\mathbb C} 
\to H^{1}({\mathbb P}^{2}_{q}, {\mathcal M}(-1)) 
\to H^{1}({\mathbb P}^{2}_{q}, {\mathcal M}) 
\to 
0\ ,
$$
whence the last statement in the Theorem.
\end{proof}
Now, as in the Introduction, let $\, V(\bs{w}) := 
H^1(\P,\, \ms{M}(-1)) \,$. It follows 
from Theorem~\ref{PP1} that multiplication by $\, Z\,$
defines {\it isomorphisms}
$$
H^1(\P,\, \ms{M}(-\bs{w})) \cong H^1(\P,\, \ms{M}(-\bs{w}+1)) 
\cong \, \ldots \, \cong H^1(\P,\, \ms{M}(-2)) \cong V(\bs{w})\ .
$$
We identify these spaces, and let $\, \X \,$ and  $\, \Y \,$
be the endomorphisms of $\, V(\bs{w})\,$ induced by (right)
multiplication by $ X $ and $ Y \,$. More precisely, if 
$\, v \in V(\bs{w}) \,$, we define 
$$
\X(v) := v \cdot Z^{-w_1} X  \quad , \quad \Y(v) := v \cdot  Z^{-w_2} Y\ .
$$
Let $\, n := \dim_{\c}V(\bs{w})\,$.  By Theorem~\ref{PP1}, we have 
$\,n = 0 \,$ if and only if $\,M\,$ is cyclic.
\begin{proposition}
\la{XY}
The pair $\, (\X, \Y) \,$ defines a point in the space 
$\, \mathfrak{C}_n \,$.
\end{proposition}
\begin{proof}
The Proposition is trivial if $\,n = 0\,$.  In general, we calculate:
$$
\X\Y(v) \cdot Z = v \cdot  Z^{-w_2} Y  Z^{-w_1} X Z = 
 v \cdot Z^{-|\bs{w}| + 1} Y X
$$
(we used the facts that $\, Z \,$ commutes with $\, X\,$ and
$\, Y\,$, and that $\, Z^{-|\bs{w}| + 1} \,$ is still well defined on 
$\,V(\bs{w})\,$). Similarly, $\, \Y\X(v) \cdot Z = v 
\cdot Z^{-|\bs{w}| + 1} XY \,$. So
$$
\left([\X,\Y] + \I\right) v \cdot Z = v \cdot Z^{-|\bs{w}| + 1} 
(YX - XY + Z^{|\bs{w}|})= 0 \ .
$$
Thus the image of $\, [\X,\Y] + \I \,$ is contained in the
kernel of $\, \cdot Z\,: V(\bs{w}) \to H^1(\P,\,\ms{M}) \,$.
By Theorem~\ref{PP1}, this map is surjective with one-dimensional
kernel. Therefore $\, [\X,\Y] + \I \,$ has rank $1\,$, as required.
\end{proof}

\section{Elementary Constructions}
\la{Sect5}

\subsection{Distinguished Representatives}
\la{DR}

As usual, let $ M $ be a finitely generated rank one 
torsion-free right $A$-module. We are going to construct 
two distinguished realizations of $ M $ as fractional ideals
of $\, A\,$ (that is, submodules of the quotient field $\, Q \,$ 
of $\, A \,$). First, according to \cite{Sta}, Lemma~4.2, we can choose
an embedding of $ M $ as an ideal which has nonzero
intersection with $\, \c[x] \subset A \,$. If an element of $ A $
(or, later, of the larger algebra $ \c(x)[y] $) is written in the  
form $\, a = \sum_{i=0}^{n} a_{i}(x) y^{i} \,$ 
(with $\, a_n \not= 0\,$), we call $\, a_{n}(x) \,$ the 
{\it leading coefficient} of $ a\,$. 
The leading coefficients of all the elements of $ M $ form
an ideal in $\, \c[x] \,$; let $\, p(x) \,$ be the (monic) 
generator of this ideal, and set $\, M_{x} := p(x)^{-1}M \,$.
By construction, the fractional ideal $\, M_{x} \,$ has the 
following properties:
\begin{enumerate}
\item $\, M_{x} \subset \c(x)[y] \subset Q \ $ and 
$\ M_{x} \cap \c[x] \not= \{0\} \,$;
\item \ all leading coefficients of elements of $\, M_{x} \,$
belong to $\, \c[x] \,$;
\item $\, M_{x} \,$ contains an element with constant 
leading coefficient.
\end{enumerate}
It is easy to see that these properties characterize $ M_x \,$.
More precisely, we have
\begin{lemma}
\la{Mx}
Let $\, M_x \,$ and $\, M_x' \,$ be two fractional ideals
of $ A \,$, both isomorphic to $ M \,$, and satisfying (1)-(3)
above. Let $\, q \,$ be an element of $\, Q \,$ such that
$\, M_x' = q M_x \,$. Then $\, q \,$ is a constant 
(and hence $\, M_x = M_x' \,$). 
\end{lemma}

We denote by $\, \varrho_{x}\,:\, \c(x)[y] \to A\, $ the map
that deletes the ``negative part'' of the coefficients 
$ a_{i}(x) \,$. More precisely, if $\, a(x)\,$ is rational function 
of $ x \,$, let $\, a = a_{+} + a_{-} \,$, where $\, a_{+} \,$
is a polynomial and $\, a_{-} \,$ vanishes at infinity; then we define
$$
\varrho_{x}\left(\sum a_{i}(x) y^{i}\right) :=  
\sum a_{i}(x)_{+}\, y^{i}\ .
$$
We denote by $\, r_{x} \,$ the restriction of $\, \varrho_{x}\, $
to $\, M_{x} \,$. Then $\, r_{x} \,$ is injective, and 
$\,r_{x}(M_{x}) \,$ is a linear subspace of finite codimension
in $ A\,$. Let $\, V_x := A/r_{x}(M_{x}) \,$. The map $\,\varrho_{x} \,$
commutes with right multiplication by $\, y \,$ (though not with 
right multiplication by $ x $); thus $\, \cdot y \,$ induces an 
endomorphism of $\, V_{x} \,$. We denote it by $\, \Y \,$.

Reversing the roles of $\, x \,$ and $\, y \,$ in the above 
construction, we obtain another distinguished representative
$\, M_{y} \subset \c(y)[x] \subset Q\,$ for our ideal $ M\,$,
and another finite-dimensional vector space
$\, V_y := A/r_{y}(M_{y})\,$, together with an endomorphism
$\, \X \,$ of $\, V_y \,$ coming from right multiplication
by $ x\,$. Since $\, M_x \,$ and $\, M_y \,$ are both isomorphic
to $ M \,$, we have $\, M_y = \kappa M_x \,$ for some 
$\, \kappa \in Q\,$; by Lemma~\ref{Mx}, $\, \kappa \,$ is uniquely
determined up to a constant factor. Note that the properties 
of $\, M_x \,$ and $\, M_y \,$ imply
\begin{equation}
\la{cxy}
\kappa \in \c(y)(x) \quad \mbox{and} \quad \kappa^{-1} \in  \c(x)(y)\ .
\end{equation}
Here $\, \c(x)(y)  \,$ (for example) denotes the space of 
all elements of $\, Q \,$ that have the form 
$\, \sum f_{i}(x)\, g_{i}(y) \,$ for some rational functions
$\, f_{i}(x) \, , \, g_{i}(y) \,$. 

Next, we describe the linear isomorphism 
$\, \phi \,:\, V_x \to V_y \,$ mentioned in the
Introduction. In the next section we shall 
see how this isomorphism arises naturally from a 
calculation of \v{C}ech cohomology.
Let $\, \Phi: r_{x}(M_x) \to r_y(M_y) \,$ be the isomorphism
defined by 
$$
\Phi(m) := r_y\left(\,\kappa \cdot r_{x}^{-1}(m)\,\right)
$$
(the dot denotes multiplication in $\,Q$). We shall extend
$\,\Phi\,$ to a linear isomorphism (also denoted by $\,\Phi$)
from $\, A \,$ to itself: $\, \phi \,$ will then be the 
induced map on the quotient spaces. The extension of 
$\, \Phi \,$ to $\, A\,$ is defined as follows.
Note first that for any $\, a\in A \,$ there are 
polynomials $\, g(y)\,$ such that $\, a \, g(y) \in r_{x}(M_x)\,$
(for example, we can take $\, g\,$ to be the characteristic 
polynomial of the map $\, \Y \,$ above). For each 
$\, a \in A \,$, we  choose such a polynomial $\, g(y)\,$,
and set
\begin{equation}
\la{map}
\Phi(a) := \varrho_{y}\left(\,\kappa \cdot r_{x}^{-1}[\,a\,
g(y)\,]\cdot g(y)^{-1}\,\right)\ .
\end{equation}
Using the $\,\c[y]$-linearity of the map $\,\varrho_{x} \,$,
it is easy to check that $\, \Phi(a) \,$ is independent of 
the choice of $\, g\,$; in particular, if $\, a \in r_{x}(M_x)\,$ 
we can choose $\, g = 1 \,$, so $\, \Phi \,$ is indeed an 
extension of the map that we started with.
The reader may like to prove at this point that 
$\,\Phi\,$ and $\, \phi\,$ are isomorphisms (with inverses
defined in a similar way, interchanging $\,x\,$ and $\,y$).
This will follow from the results of the next section, so 
we omit the proof here. 

\subsection{An Example}
\la{EX}
For $\, n \geq 1 \,$, let $\, M = x^{n+1} A + (xy+n) A\,$. 
In this case the spaces $\, V_x \,$ and $\, V_y \,$ have dimension
$\, n \,$, and (with suitable choice of basis) the
matrices $\, \X \,$ and $\, \Y \,$ are
$$
\X \,=\, \left(
\begin{array}{ccccc}
0 & 0 &  0 & \ldots & 0\\
1 & 0 &  0 & \ldots & 0\\
0 & 1 &  0 & \ddots & \vdots\\
\vdots & \vdots & \ddots & \ddots & 0 \\
0 & 0 & \ldots & 1 & 0
\end{array}
\right)\quad , \quad
\Y \, = \, \left(
\begin{array}{ccccc}
0 &  1-n &     0  & \ldots & 0\\
0 &    0  &   2-n & \ldots & 0\\
0 &    0  &    0  & \ddots & \vdots\\
\vdots & \vdots & \ddots & \ddots & -1 \\
0 & 0  & \ldots & 0 & 0
\end{array}
\right)\ .
$$
Although very elementary, the calculation is not short enough 
to reproduce in full here. We just indicate the main steps,
leaving some details for the reader.
First, we have $\, M_x = x^{-1} M \,$; as a basis for
$\, V_x \,$ we can take the residue classes 
(modulo $\, r_x(M_x)$) of the elements $\, 1, x,\ldots, x^{n-1} \,$. 
Since $\, x^{k}y + (n-k) x^{k-1} \, \in \, r_x(M_x)\,$ 
for $\, 1 \leq k \leq  n-1 \,$, and $\,y \in r_x(M_x)\,$, 
it follows that the matrix 
$\, \Y \,$ is as above. The calculation of $\, \X \,$ is 
a little harder; however, using formula (\ref{map}), it is
straightforward to check that
$\, \Phi(x^k) = x^k \,$ for  $\, 0 \leq k \leq  n-1 \,$.
So we can again choose the residue classes of 
$\,1, x,\ldots, x^{n-1} \,$ (now modulo $\, r_y(M_y)$) 
as a basis for $\, V_y \,$, and the matrix of 
$\, \phi \,$ is then the identity. Since 
$\, x^n \in r_y(M_y)\,$, it follows that the matrix 
$\, \X \,$ is as above. In this example we have
$$
\kappa = (xy)^{-n} y\,(xy+1)\,(xy+2)\,\ldots \, (xy+n-1)\,x
\ \ \mbox{and}\ \ 
\kappa^{-1} = 1 + n\,x^{-1}y^{-1}\ .
$$

\subsection{The Associated Graded Ideals}
\la{AGI}

Our last goal in this section is to establish an important
property of the element $\,\kappa \,$ which we shall need later,
namely, that multiplication by $\,\kappa \,$ preserves the 
$ \bs{w} $-filtration on $ \,A\, $ for every weight vector 
$ \bs{w} = (w_1 , w_2) \, $.
Slightly more generally than in 
Section~\ref{SBf}, we shall allow $\, w_1 \,$ and $\, w_2 \,$ 
to be any non-negative integers\footnote{That is, we now allow
one of $\, w_{i} \,$ to be zero. More generally still, we could 
work with non-integer (real) ``weights'' as in \cite{Di}, \cite{LM}.} 
that are not both zero.
We denote by $\, \vv_{\bs{w}} \,$ the valuation on $\, A \,$
corresponding to $\,\bs{w}\,$, that is, if $\, a \not=0 \,$ 
then $\, \vv_{\bs{w}}(a) \,$ is the least integer $\, k \,$ 
such that $\, a \in A_{k}(\bs{w}) \,$ (and 
$\, \vv_{\bs{w}}(0) := -\infty \,$). We extend $\, \vv_{\bs{w}} \,$
to $\, Q\,$ by setting
$$
\mathfrak{v}_{\boldsymbol{w}}(ab^{-1}) = 
\mathfrak{v}_{\boldsymbol{w}}(a) - \mathfrak{v}_{\boldsymbol{w}}(b)\ ,
$$
and let $\, Q_{k}(\bs{w}) := \{ \, 
q \in Q\ |\ \mathfrak{v}_{\boldsymbol{w}}(q) \leq k\,\} \,$.
Then $\, \{Q_{\bullet}(\bs{w})\} \,$ is a separated 
filtration on $ \, Q \,$ extending the original filtration
on $\, A\,$. Changing notation slightly from Section~\ref{SBf},
we denote the associated graded algebra by $\,\bs{G_{w}Q}\,$
and write
$$
\sigma_{\bs{w}}\,:\, Q \to \bs{G_{w}Q}
$$
for the {\it symbol map}\,:\, if $\, \vv_{\bs{w}}(q) = k \,$, 
then $\,\sigma_{\bs{w}}(q) \,$ is the class of $\, q \,$ in 
$\,Q_{k}(\bs{w})/Q_{k-1}(\bs{w})\,$. As usual, we identify 
$\, \bs{G_{w}A} \equiv \bs{GA}(\bs{w}) \,$ with the polynomial algebra
$\, \c[\bar{x}, \bar{y}]\, $, where $\,\bar{x} :=  \sigma_{\bs{w}}(x)\,$ 
and $\,\bar{y} := \sigma_{\bs{w}}(y)\,$; then
$\, \bs{G_{w}Q} \,$ is identified with the subalgebra of 
$\, \c(\bar{x}, \bar{y})\,$ spanned by quotients of homogeneous 
polynomials. 

For short, we write $\,\sigma_{y} \,$ and $\,\bs{G}_{y} \,$
instead of $\, \sigma_{(0,1)}\,$ and $\,\bs{G}_{(0,1)}\,$.
Property $ (2) $ of $\, M_x = p(x)^{-1}M \,$ says that
$\, \bs{G}_{y}M_x \subseteq \c[\bar{x}, \bar{y}] \,$
(even though $\, M_x \not\subset A $). More precisely,
if $\, a_{n}(x) \,$ is the leading coefficient of 
$\, a \in M \,$ then $\, \sigma_{y}(a) = a_{n}(\bar{x})\bar{y}^n \,$.
It follows that $\, p(\bar{x}) \,$ is the greatest common divisor
of the elements of $\, \bs{G}_{y}M\,$ so that 
$\, \bs{G}_{y}M_x =  p(\bar{x})^{-1} \bs{G}_{y}M \,$ is an ideal
of finite codimension\footnote{equal to the codimension of
$\, r_x(M_x)\,$ in $\,A\,$.} in $\, \c[\bar{x}, \bar{y}]\,$.
More generally, we have
\begin{proposition}
\la{Mak}
Let $\,\bs{w}\,$ be any weight vector (as specified above). 
Then

$ (i)$\ $\, \bs{G_w}M_x \,$ and $\, \bs{G_w}M_y \,$ are ideals 
of finite codimension in $\,\c[\bar{x}, \bar{y}]\,$;

$ (ii) $\ $\, \bs{G_w}M_x =  \bs{G_w}M_y\,$ in $\,\c[\bar{x}, \bar{y}]\,$;

$ (iii) $\ the symbol $\, \sigma_{\bs{w}}(\kappa) \,$ is constant.
\end{proposition}
\begin{proof}
An argument of Letzter and Makar-Limanov (see \cite{LM}, Lemma~2.1) 
shows that $\, \bs{G_w}M_x \,$  is contained in 
$\, \c[\bar{x}, \bar{y}]\,$. Proposition~2.4$'$ of \cite{LM} 
then shows that it has finite codimension (and also
that this codimension is independent of $\bs{w}\,$).
Interchanging the roles of $\, x \,$ and  $\, y \,$, 
we obtain the same result for $\, \bs{G_w}M_y\,$.

Since $\, M_x \,$ and $\, M_y \,$ are both isomorphic to $ M \,$,
the ideals $\,\bs{G_w}M_x\,$ and $\, \bs{G_w}M_y\, $ are isomorphic.
An ideal class of $\, \c[\bar{x}, \bar{y}] \,$ has a {\it unique} 
representative of finite codimension, so $(ii)$ follows from $(i)$.

Since $\, M_y = \kappa M_x \,$,
we have $\, \bs{G_{w}} M_y = \sigma_{\bs{w}}(\kappa) \bs{G_{w}}M_x \,$
(the symbol map is multiplicative). So $ (iii) $ follows from $(ii)$.
\end{proof}
\begin{corollary}
\la{anyind}
For any positive weight vector $\,\bs{w}\,$, the filtration induced 
on $\,M_x\,$ (or $\,M_y\,$) by the $\bs{w}$-filtration on $\,Q\,$ 
coincides with the {\it normalized} induced filtration of 
Section~\ref{Sect4}.
\end{corollary}
\begin{proof}
This is a reformulation of Proposition~\ref{Mak}$ (i) $ (cf.\ the proof of 
Lemma~\ref{LL1} above).
\end{proof}
\subsection*{Remark} 
If we compare the formula $\, M_y = \kappa M_x \,$ with 
Proposition~6.2 in \cite{BW}, we see that $ \,\kappa \,$
can be identified with the formal integral operator
$\, K \,$ that plays a basic role in the theory of
integrable systems; more precisely, if $\, W \,$ is
the point of $\, \mbox{\rm Gr}^{\mbox{\rm \scriptsize ad}} \,$
that corresponds to the ideal $\, M \,$ then 
$\, \kappa = K_{b(W)} \,$, where $\, b \,$ is the 
bispectral involution on 
$\, \mbox{\rm Gr}^{\mbox{\rm \scriptsize ad}} \,$.
We shall not make any use of this remark in the present
paper; however, it points the way to a more direct proof
of Theorem~\ref{IT2}.

\section{The Comparison Theorem}
\la{Sect6}

As usual, let $\, \ms{M}\,$ be the canonical extension of a
noncyclic ideal $ M $ of $ A\,$. Our aim in this section is to 
calculate the groups $\, H^1(\P(\bs{w}), \, \ms{M}(k)) \,$
(for $\, -|\bs{w}| \leq k \leq -1 $) using the 
\v{C}ech complex of the covering $\, \U \,$ introduced in 
Section~\ref{SBf}: this will enable us to
identify these groups with the spaces $\, V_x \,$
and $\, V_y \,$ in Section~\ref{Sect5}. Although it would, 
of course, be possible to calculate 
$\, \H^1(\U, \, \ms{M}) \,$ directly from the complex 
(\ref{10}), this does not appear to yield the answer in 
the form we want. Instead, we choose a large integer 
$\, p \,$ (eventually we shall let $\, p \to \infty$), and
denote by $\, \ms{N}_p \,$ the restriction of 
$\, \ms{M} \,$ to the $ p$-th infinitesimal neighbourhood
of the line at infinity in $\, \P\,$; that is, 
$\, \ms{N}_p\,$ is the quotient term in the exact sequence
\begin{equation}
\la{13} 
0 \to \ms{M} \stackrel{\cdot Z^p}{\longrightarrow} \ms{M}(p) \to 
\ms{N}_p \to 0\ .
\end{equation}
In what follows we shall mostly omit the subscript $\, p \,$,
writing $\, \ms{N} \equiv \ms{N}_p \,$. 
We shall assume that $\, p \,$ is chosen so that
\begin{equation}
\la{133}
H^{1}(\P,\, \ms{M}(p+k)) = 0 
\quad \mbox{for all}\ k \geq - |\bs{w}|
\end{equation}
(this is possible by the Vanishing Theorem~\ref{FVT}$(b)$).
From (\ref{13}), we then get the exact sequence 
\begin{equation}
\la{15} 
0 \to  H^{0}(\P,\, {\mathcal M}(p+k)) 
  \to  H^{0}(\P,\, {\mathcal N}(k)) 
  \to  H^{1}(\P,\, {\mathcal M}(k)) 
  \to  0 
\end{equation}
for any $\, k \,$ in the range $\, -|\bs{w}| \leq k \leq -1 \,$. 
We are going to calculate
$\, H^{0}(\P,\, {\mathcal N}(k)) \,$ via a \v{C}ech complex,
and then obtain $\, H^{1}(\P,\, {\mathcal M}(k)) \,$
as the quotient term in (\ref{15}). 
We first record the following fact.
\begin{lemma}
\la{DIM}
For any $\, k \geq - |\bs{w}| \,$, we have
$$\
\dim_{\c}  H^{0}(\P,\, {\mathcal N}(k)) = \dim_{\c} 
A_{p+k} - \dim_{\c} A_k\ .
$$
\end{lemma}
\begin{proof} 
From (\ref{13}) we easily find (using (\ref{133}) and Theorem~\ref{PP1})
that 
$$ 
H^{1}(\P,\, \ms{N}(k)) =  H^{2}(\P,\, \ms{N}(k)) = 0 
\quad \mbox{for all}\ k \geq  -|\bs{w}|\ ,
$$
so $\, \dim_{\c}  H^{0}(\P,\, \ms{N}(k)) = \chi(\P,\,\ms{N}(k))\,$
(where $\, \chi \,$ denotes the Euler characteristic). 
From (\ref{13}) again, we then get
\begin{equation}
\la{DIM0}
\dim_{\c}  H^{0}(\P,\, {\mathcal N}(k)) = \chi(\P,\,\ms{M}(p+k)) - 
\chi(\P,\,\ms{M}(k))\ .
\end{equation}
On the other hand, by (\ref{5.3}) the Euler characteristics of 
the sheaves $\, \ms{M}(r)\,$
satisfy
$$
\chi(\P,\,\ms{M}(r)) = \chi(\P, \, \ms{M}(r-1)) + 
\chi(\PP^{1},\,\ms{O}_{\PP^{1}}(r))\ ;
$$
and by (\ref{5.14}), $\, \chi(\PP^{1},\,\ms{O}_{\PP^{1}}(r)) = 
\dim_{\c} S_{r} = \dim_{\c} A_{r} - 
\dim_{\c} A_{r-1}\,$ for all $\, r > -|\bs{w}| \,$.
By Theorem~\ref{PP1}, we have $\,\chi(\P,\,\ms{M}(-|\bs{w}|)) 
= -n \,$, so by induction 
\begin{equation}
\la{DIM2}
\chi(\P,\, {\mathcal M}(r)) = \dim_{\c} A_{r} - n \quad \mbox{for all}
\hspace*{0.2cm} r \geq -|\bs{w}| \ .
\end{equation}
Combining (\ref{DIM0}) and (\ref{DIM2}) yields the Lemma.
\end{proof}

For a while now  we shall work with the special representative
$\, M_x \,$ for the class of $\, M\,$ (see Section~\ref{Sect5});
to simplify the notation we drop the suffix $\, x \,$ and 
denote $\, M_x \,$ simply by $\, M\,$. We have the (normalized)
filtration on $\, M\,$ induced by the $ \bs{w}$-filtration on 
$\, Q\,$ (see Corollary~\ref{anyind}). 
As usual, let $\, \bs{M} = \bigoplus_{k \in \Z} M_k \,$
be the corresponding homogenization of $\, M\,$, and let $\, \bs{N}\,$ 
denote the quotient term in the exact sequence
\begin{equation}
\la{21}
0 \to \bs{M} \stackrel{\cdot Z^p}{\longrightarrow} 
\bs{M}(p) \to \bs{N} \to 0 \quad  ,
\end{equation}
so that $\, \pi \bs{N} = \ms{N}\,$. Let $\, \mathfrak{U} = 
(U_1, U_2)\,$ be the covering of $\, \bs{A}\,$ defined by
(\ref{4}). We identify the localizations 
$\, \bs{A}_{U_1}\,$ and $\, \bs{A}_{U_2}\,$ with subalgebras
of the homogenization $\, \bs{Q} \,$ of the Weyl quotient 
field $\, Q\,$: specifically, we have 
$$
\bs{A}_{U_{1}} = \bigoplus_{k \in \mathbb{Z}} 
\mathbb{C}(x)[y]_k \,  \quad \mbox{and} \quad
\bs{A}_{U_{2}} = \bigoplus_{k \in \mathbb{Z}} 
\mathbb{C}(y)[x]_k \,  \ 
$$
(here and below, the subscript $ k $ refers to the filtration induced 
from $\, Q \,$).    
In a similar way, the embedding of $\, M \,$ in $\, Q\,$
allows us to identify the localizations of $\, \bs{M} \,$
with subspaces of $\,\bs{Q}\,$ (the tensor products in 
Definition~\ref{8} then become multiplication in $\,\bs{Q}$).
First, the conditions $\, M \cap \c[x] \not= \{0\} \,$,
$\, M \subset \c(x)[y] \,$ imply that 
$\, \bs{M}_{U_1} = \bs{A}_{U_1} \,$. To calculate 
$\, \bs{M}_{U_2} \,$, we use the other distinguished
representative $\, M_y = \kappa\, M \,$.
For the corresponding homogenizations, we have
$\, \bs{M}_y = \bs{\kappa}\,\bs{M} \,$.
As above, the localization of $\, \bs{M}_y \,$ with 
respect to $\, U_2 \,$ is just $\, \bs{A}_{U_2} \,$;
and by Proposition~\ref{Mak}$(iii)$, multiplication by 
$\, \kappa \,$ preserves the $ \bs{w}$-filtration 
on $\, Q \,$, that is $\, \bs{\kappa} \in Q_0\,$. 
We therefore have 
$$
\bs{M}_{U_{1}} = \bigoplus_{k \in \mathbb{Z}} 
\mathbb{C}(x)[y]_k \,  \quad \mbox{and} \quad
\bs{M}_{U_{2}} = \bigoplus_{k \in \mathbb{Z}} 
\kappa^{-1} \mathbb{C}(y)[x]_k \,  \ .
$$
Finally, because localization is an exact functor,
we can calculate $\, \bs{N}_{U_{1}} \,$ and 
$\, \bs{N}_{U_{2}} \,$ by localizing the exact 
sequence (\ref{21}). The result is
$$
\bs{N}_{U_{1}} = \bigoplus_{k \in \mathbb{Z}} 
\, \mathbb{C}(x)[y]_{p+k}/\mathbb{C}(x)[y]_k  
\ , \
\bs{N}_{U_{2}} = \bigoplus_{k \in \mathbb{Z}} 
\, \kappa^{-1} \mathbb{C}(y)[x]_{p+k}/\kappa^{-1} 
\mathbb{C}(y)[x]_k  \ .
$$
The repeated localizations $\, \bs{N}_{U_i U_j} \,$ similarly
get identified with certain (easily specified) subspaces of
$\, \bigoplus_{k \in \Z} Q_{p+k}/Q_{k}\,$. In what follows we 
shall denote elements of degree $\, k \,$ in 
$\, \bs{N}_{U_{1}} \,$ and $\,\bs{N}_{U_{2}}\,$ 
by $\, \bar{n}_1 \,$ and $\, \bar{n}_2 \,$, where it is understood that 
$\, n_1 \in \mathbb{C}(x)[y]_{p+k}\,$,
$\, n_2 \in \kappa^{-1} \mathbb{C}(y)[x]_{p+k}\,$, and that 
the bars denote residue classes modulo elements in $Q_k\,$.
\begin{proposition}
\la{hn}
With the identifications explained above, 
$\, \check{H}^{0}(\mathfrak{U}, \ms{N}(k)) \,$
is the subspace of $\,(\bs{N}_{U_{1}})_k \oplus (\bs{N}_{U_{2}})_k
\,$ consisting of all pairs $\, (\bar{n}_1, \bar{n}_2)\, $ such that 
$\, n_1- n_2 \in Q_{k}\,$. Furthermore, the map
$\, M_{p+k} \cong \check{H}^{0}(\mathfrak{U}, \ms{M}(p+k)) 
\to \check{H}^{0}(\mathfrak{U}, \ms{N}(k)) \,$
coming from (\ref{13}) sends $\, m \,$ to 
$\, (\bar{m}, \bar{m}) \,$.
\end{proposition}
\begin{proof}
The first statement is obvious, because the coboundary map
$$
d^{\,0}\,:\ \bs{N}_{U_{1}} \oplus \bs{N}_{U_{2}} \to 
\bs{N}_{U_{1}U_{1}} \oplus \bs{N}_{U_{1}U_{2}} \oplus \bs{N}_{U_{2}U_{1}} 
\oplus \bs{N}_{U_{2}U_{2}}
$$
takes $\, (\bar{n}_1, \bar{n}_2) \,$ to 
$\, (0, \bar{n}_1-\bar{n}_2,  \bar{n}_2-\bar{n}_1, 0) \,$. 
The second statement follows from Proposition~\ref{ch}.
\end{proof}

Now, for each $\, k \in \Z\,$, define a map
\begin{equation}
\la{gk}
\gamma_k :\, \check{H}^{0}(\mathfrak{U}, \ms{N}(k)) \to
A_{p+k}/A_{k}
\end{equation}
by setting $\, \gamma_k (\bar{n}_1, \bar{n}_2) := 
\overline{\varrho_x(n_1)} \,$, where 
$\,\varrho_x :\, \c(x)[y] \to A \,$ is as in Section~\ref{Sect5}.
\begin{proposition}
\la{isogk}
The map $\, \gamma_k \,$ is an isomorphism for all 
$\, k \geq - |\bs{w}|\,$.
\end{proposition}
\begin{proof}
By Lemma~\ref{DIM}, the two spaces in (\ref{gk}) have the same
(finite) dimension if $\, k \geq - |\bs{w}|\,$, so it is enough
to prove that $\, \gamma_{k} \,$ is injective if 
$\, k \geq - |\bs{w}|\,$. In fact, $\, \gamma_k \,$ is injective
for all $ k \in \Z\,$. To see that, let
$\, (\bar{n}_1, \bar{n}_2) \in 
\check{H}^{0}(\mathfrak{U}, \ms{N}(k)) \,$ and suppose
$\, \gamma_k (\bar{n}_1, \bar{n}_2) = 0 \,$; we have
to show that $\, \bar{n}_1 = \bar{n}_2 = 0 \,$.
Equivalently, by Proposition~\ref{hn}, we are 
given $\, n_1 \in \c(x)[y]_{p+k} \,$
and $\, n_2 \in \kappa^{-1}\c(y)[x]_{p+k} \,$ such that
$\, n_1 - n_2 \in Q_k \,$ and $\, \varrho_x(n_1) \in Q_k \,$; 
we have to show that $\, n_1 \,$ and $\, n_2 \,$  are in 
$\, Q_k\,$. Clearly, it is enough to show that 
$ n_2 \in  Q_k\,$. We extend $\, \varrho_x \,$ to a 
map from $\,\c(x)(y)\,$ to $\,\c[x](y) = \c(y)[x] \,$ by setting
(as in Section~\ref{Sect5})
$$
\varrho_x \left(\sum f_i(x) g_{i}(y)\right) := 
\sum f_i(x)_{+}\, g_{i}(y)\ .
$$
It is easy to see that $\, \varrho_x \,$ is well defined
and respects the $ \bs{w}$-filtration for any $\,\bs{w}\,$.
Note that $\, \kappa^{-1} \c(y)[x] \subset \c(x)(y) \,$ by 
(\ref{cxy}), so $\, n_2 \in \c(x)(y) \,$. We have
$\, n_1 - n_2 \in Q_k \,$, hence 
$\, \varrho_x(n_1) - \varrho_x(n_2) \in Q_k \,$; 
since we are given $\, \varrho_x(n_1) \in Q_k \,$,
we get $\, \varrho_x(n_2) \in Q_k \,$. But we claim
that if $\, n \in \kappa^{-1} \c(y)[x] \,$ then 
$\, \vv_{\bs{w}}(n) = \vv_{\bs{w}}(\varrho_x(n)) \,$,
hence $\, n_2 \in Q_k\,$, as required.
To prove the last claim, write $\, n = \kappa^{-1} q \,$ with $\, q \in 
\c(y)[x] \,$. By Proposition~\ref{Mak}$(iii)$, we may normalize
$\, \kappa \,$ so that $\, \sigma_{\bs{w}}(\kappa) = 1\,$.
We then have $\, n = q + q' \,$, where 
$\, \vv_{\bs{w}}(q') < \vv_{\bs{w}}(q)\,$;
hence $\, \vv_{\bs{w}}(n) = \vv_{\bs{w}}(q)\,$. Since
$\, q \in \c(y)[x]\, \Rightarrow \,\varrho_x(q) = q \,$, we have 
$\, \varrho_x(n) = q \, + \, \varrho_x(q')\,$, and 
$\, \vv_{\bs{w}}(\varrho_x(q')) \leq  \vv_{\bs{w}}(q') < 
\vv_{\bs{w}}(q) \,$. Hence $\, \vv_{\bs{w}}(\varrho_x(n)) = 
\vv_{\bs{w}}(q) = \vv_{\bs{w}}(n)\,$,
as claimed above. 
\end{proof}
Combining the maps $\, \gamma_k \,$ for all $\, k \in \Z \,$,
we now get an (injective) map
\begin{equation}
\la{grga}
\bs{\gamma}:\, \underline{\check{H}}{}^{0}(\mathfrak{U}, \ms{N}_p) 
\to \bigoplus_{k\in \Z} A_{p+k}/A_{k}\ .
\end{equation}
The $\, \c[y]$-linearity of the map $\,\varrho_x \,$
implies that $\, \bs{\gamma} \,$ is a homomorphism
of graded $\, \c[Y,Z]$-modules, where on the right
$\, Y \,$ acts as right multiplication by $\, y \,$
and $\, Z \,$ acts by embedding successive filtration
components.

We now focus our attention on the degrees $\, k \,$
in the range $\, -|\bs{w}| \leq k \leq -1 \,$. In this 
case $\, A_k = 0 \,$, so $\, \gamma_k \,$ is simply an 
isomorphism $\, 
\check{H}^{0}(\mathfrak{U}, \ms{N}(k)) 
\to A_{p+k} \,$. Looking back at the exact sequence 
(\ref{15}) and using the last statement in 
Proposition~\ref{hn}, we see that $\, \gamma_k \,$
induces isomorphisms
\begin{equation}
\la{hone}
H^{1}(\P,\, \ms{M}(k)) \stackrel{\sim}{\longrightarrow} 
A_{p+k}/r_{x}(M_{p+k}) \quad \mbox{for}\  -|\bs{w}| \leq k \leq -1
\end{equation}
(recall that $\, r_{x} \,$ denotes the restriction of $\, \varrho_x \,$ 
to $\, M$). 
We can now let $\, p \to \infty \,$. It is easy to check that 
the map $\, 
\underline{\check{H}}{}^{0}(\mathfrak{U}, \ms{N}_{p}) \to 
\underline{\check{H}}{}^{0}(\mathfrak{U}, \ms{N}_{p+1})\,$
induced by multiplication by $\, Z \,$ is compatible with
embedding of components $\, A_{p+k}/A_k \hookrightarrow 
A_{p+k+1}/A_k\,$ on the right of (\ref{grga}). It follows 
that the isomorphisms (\ref{hone}) are compatible with
the embeddings $\, A_{p+k}/r_{x}(M_{p+k}) \hookrightarrow 
A_{p+k+1}/r_{x}(M_{p+k+1})\,$; hence, letting 
$\, p \to \infty \,$ in (\ref{hone}), we get isomorphisms
\begin{equation}
\la{alp}
\alpha_k :\, 
H^{1}(\P,\, \ms{M}(k)) \to V_x\ , \quad  -|\bs{w}| \leq k \leq -1\ ,
\end{equation}
where (as in Section~\ref{Sect5}) $\, V_x := A/r_{x}(M)\,$.
Further, the $\, \c[Y,Z]$-linearity of $\, \bs{\gamma} \,$
implies that the isomorphisms $\, \alpha_k \,$ take 
multiplication by $ Y $ and $ Z $ (when defined) on the left
of (\ref{alp}) to (right) multiplication by $\, y \,$ and to the 
identity map (respectively) on the right. It follows at once
that the isomorphism 
$$
\alpha_x :=  \alpha_{-1}:\, H^{1}(\P,\, \ms{M}(-1)) = 
V(\bs{w}) \to V_x
$$
takes the map $\, \Y(\bs{w})\,$ of Section~\ref{Sect4} to
the map $\, \Y \,$ of Section~\ref{Sect5}, as claimed in 
Theorem~\ref{IT5}. 

To obtain the other isomorphism $\, \alpha_y \,$ in 
Theorem~\ref{IT5}, we have only to repeat all the above,
starting from the representative $\, M_y \,$ rather than
$\, M_x \,$. We sketch a few details to fix the notation
for the last calculation below. To avoid confusion, we 
continue to denote $\, M_x \,$ by $\, M \,$, so that
$\, M_y = \kappa\,M\,$. Further, we continue to identify
$\, \check{H}^{0}(\mathfrak{U}, \ms{N}(k)) \,$
with the space described in Proposition~\ref{hn}, so that
in the new argument we work with the realization 
$\, \kappa \,\check{H}^{0}(\mathfrak{U}, \, \ms{N}(k)) \,$.
The crucial map
$$
\gamma_{k}' :\, \kappa\,
\check{H}^{0}(\mathfrak{U}, \, \ms{N}(k)) \to A_{p+k}/A_{k}
$$
is then defined by  $\, \gamma_{k}'(\kappa\,\bar{n}_1, 
\kappa\,\bar{n}_2) := \overline{\varrho_y(\kappa\,n_2)} \,$.
Passing to a quotient and letting $\, p \to \infty \,$, we get
the required isomorphism $\, \alpha_y : V(\bs{w}) \to V_y \,$
exactly as before.

To complete the proof of Theorem~\ref{IT5}, it remains to show
that the isomorphism $\,\alpha_{y}\, \alpha_{x}^{-1}: 
\, V_x \to V_y \,$ coincides with the map $\, \phi \,$
in Section~\ref{Sect5}. To do that, we return temporarily
to the case of finite $\, p \gg 0 \,$ and 
(for $ -|\bs{w}| \leq k \leq -1$) let 
$\, \Phi_{k} \,$ be the map that makes the diagram
$$
\begin{CD}
M_{p+k} @>>>  \check{H}^{0}(\mathfrak{U}, \ms{N}(k)) @> 
\bs{\gamma}_{k} >> A_{p+k} \\
@VVV @VVV  @VV \Phi_{k} V  \\
\kappa\,M_{p+k} @>>> \kappa\,\check{H}^{0}(\mathfrak{U}, 
\ms{N}(k)) @> \bs{\gamma}_{k}' >> A_{p+k}
\end{CD}
$$
commutative. In this diagram the first two vertical arrows
are just multiplications by $\, \kappa \,$; the 
horizontal maps $\, M_{p+k} \to  A_{p+k} \,$ and 
$\, \kappa\, M_{p+k} \to  A_{p+k} \,$ are
$\, r_x \,$ and $\, r_y \,$, respectively. 
Let $\, a \in A_{p+k} \,$, and (as in Section~\ref{Sect5})
choose a polynomial $\, g(y) \,$ so that $\, a\,g(y) \in 
r_x(M) \,$, say $\, a\,g(y) = r_x(m)\,$: so here 
$\, m \in M_{p+k+N} \,$, where $\, N \,$ is the (weighted)
degree of $\, g \,$. Let $\, \gamma_{k}^{-1}(a) = 
(\bar{n}_1\,,\,\bar{n}_2)\,$; then $\,\gamma_{k}^{-1}(a\,g(y)) = 
(\overline{n_1\,g(y)}\,,\, \overline{n_2\,g(y)})\,$. On the 
other hand (by the last assertion in Proposition~\ref{hn})
$\, \gamma_{k}^{-1}(a\,g(y)) = (\bar{m}\,,\, \bar{m})\,$,
hence $\, n_2 \, g(y) - m \in Q_{k+N} \,$. Multiplying
on the left by $\, \kappa \in Q_0 \,$ and on the right
by $\, g(y)^{-1} \in Q_{-N} \,$, we get
\begin{equation}
\la{qk}
\kappa\, n_2\, - \kappa\,m \, g(y)^{-1} \in Q_{k}\ .
\end{equation}
Note that both terms in (\ref{qk}) belong to $\, \c(y)[x] \,$.
Now we can calculate:
\begin{eqnarray}
\Phi_{k}(a) 
&=& \gamma_{k}'(\kappa\,\gamma_{k}^{-1}(a)) \nonumber\\*[1ex]
&=& \gamma_{k}'\left(\kappa\,
\bar{n}_1\,,\,\kappa\,\bar{n}_2\right)\nonumber\\*[1ex]
&=& \varrho_{y}(\kappa\, n_2)  \nonumber\\*[1ex]
&=& \varrho_{y}\left(\kappa\,m \, g(y)^{-1}\right)
\hspace*{1.5cm}(\,\mbox{by}\ 
(\ref{qk})\,) \nonumber\\*[1ex]
&=& \varrho_{y}\left(\kappa\,  
r_{x}^{-1}[\,a\,g(y)\,]\, g(y)^{-1}\right)\ . \nonumber
\end{eqnarray}
Letting $\, p \to \infty \,$ we get the isomorphism 
$\, \Phi: A \to A \,$ already defined in 
Section~\ref{Sect5} (cf. formula (\ref{map})).
It follows at once that 
$\, \alpha_{y}\,\alpha_{x}^{-1} = \phi \,$, because
these are both derived from $\, \Phi \,$ by passing to 
the quotients. That completes the proof of 
Theorem~\ref{IT5}; as explained in the Introduction,
the Comparison Theorem~\ref{IT4} follows immediately.

\section{Proof of Theorem~\ref{IT3} and Theorem~\ref{IT2} }
\la{Sect7}
The natural action of $\, G = \mbox{\rm Aut}(A) \,$
on $\, \mathfrak{R} \,$ can be defined in two (equivalent) ways.  First,  
if $\, M \subseteq A \,$ is an
{\it embedded} ideal, we can make $\, \sigma \in G \,$
act on $\, M\,$ pointwise: 
$\, \sigma(M) = \{\,\sigma(m) \ | \ m \in M \,\}\,$. 
This leads to a well defined (left) action of $\, G \,$
on the space of isomorphism classes $\, \mathfrak{R} \,$, 
and is the definition used in \cite{BW}.
Alternatively, we have simply $\, \sigma(M) \cong \rho_{*}M\,$,
where $\, \rho = \sigma^{-1} \,$ is the automorphism
inverse to $\,\sigma \,$. 

Now suppose that $\sigma$  preserves the $\bs{w}$-filtration on $\, A \,$
 for some weight vector $\, \bs{w} \,$.  Then $\, \rho \,$  extends to a
graded automorphism $\, \bs{\rho} \,$ of $\,\bs{A} = \bs{A}(\bs{w})\,$, 
and the functor
$\, \bs{\rho_{*}} : \gma \to \gma \,$ descends to a functor $\, \rho_{*} \,$ 
on the quotient category $\, \C{\P(\bs{w})} \,$. 
In general, there will be no $ \bs{w}$-filtration preserved by $ \sigma $:
however, that is the case if $\sigma$ is one of the generators of $\, G \,$
in (\ref{auto}).  Slightly more generally\footnote{We need this generality
to deal with the case when $r$ or $s$ in (\ref{auto}) is 0.}, 
let $\, \sigma \,$ be an automorphism of $\, A \,$ of the form
$$
\sigma(x) = x \ , \quad \sigma(y) = y + f(x) \ ,
$$
where $\, f(x) = \sum_{0}^{r} a_i x^i \,$ is a polynomial 
in $x$ of degree $r$.  Let $\, \bs{w} = (1, N)\,$, where 
$\, N \geq r \,$.  Then $\sigma$ preserves the $\bs{w}$-filtration, 
and the extension of $\, \rho = \sigma^{-1} \,$ to $\, \bs{A} \,$ is  
defined on generators by the formulas
$$
\bs{\rho}(X) = X\ ,\quad 
\bs{\rho}(Y) = Y - \sum_{i = 0}^{r} a_i X^i Z^{N - i} \ ,\quad 
\bs{\rho}(Z) = Z\ . 
$$
If $\, M \,$ is an ideal of $\, A\, $ and  $\, \ms{M} \,$
is its canonical extension, then it is easy to see that 
the canonical extension
of $\, \sigma(M) = \rho_{*}M \,$ is $\, \rho_{*}\ms{M} \,$.
Now, for any $\, \ms{F} \in \C{\P(\bs{w})} \,$ there are natural
isomorphisms of graded $ \bs{A}$-modules
$$
\underline{H}^{i}(\P,\,\rho_{*}\ms{F}) \cong  
\bs{\rho_{*}} \underline{H}^{i}(\P,\, \ms{F}) 
$$
(cf.\ \cite{AZ}, p.~283). Thus we may identify 
$\, \underline{H}^{1}(\P,\, \rho_{*}\ms{M}) \,$ with 
$\, \underline{H}^{1}(\P,\, \ms{M})\,$ as a graded 
vector space, but multiplication by $\, X,\, Y,\, Z\,$ on 
$\,\underline{H}^{1}(\P,\, \ms{M})\,$ is then replaced by
multiplication by $\, \bs{\rho}(X),\, \bs{\rho}(Y),\, 
\bs{\rho}(Z)\,$. It follows at once that if 
$\,(\X,\,\Y)\,$ is the pair of matrices associated to 
$\, M\,$ by the construction of 
Section~\ref{Sect4}, then the pair associated to $\, \sigma(M) \,$
is $\,(\X, \, \Y - f(\X))\,$.
A similar argument (interchanging the roles of $x$ and $y$) 
shows that if $ \sigma $ is an automorphism of the form
$$
\sigma(x) = x + g(y) \ , \quad \sigma(y) = y \ ,
$$
then $ \sigma $ sends $\,(\X,\,\Y)\,$ 
to $\,(\X - g(\Y),\,\Y)\,$. These are exactly the 
formulas that defined the action of $\, G\,$ on 
$\, \mathfrak{C} \,$ in \cite{BW}, so the proof of 
Theorem~\ref{IT3} is complete.
\subsection*{Remark}
This action of $\, G\,$ on $\, \mathfrak{C} \,$ perhaps 
deserves comment, since it is not immediately obvious 
that it is well defined.  Indeed, 
if $\, \sigma \in G \,$ we are proposing to define
$\, \sigma(\X,\,\Y) \,$ by writing $\, \sigma \,$
as a product of generators $\, \Psi_{r,\lambda }\,$ and 
$\,\Phi_{s,\mu}\,$; since the matrices $\,(\X,\,\Y) \,$ 
do not satisfy the defining relation of the algebra 
$\, A \,$, it is not {\it a priori} clear that the 
result is independent of the choice of the 
representation for $\, \sigma \,$. The best way out
of this difficulty is to appeal to a theorem of 
Makar-Limanov~\cite{M} which implies that the relations
satisfied by $\, \Psi_{r,\lambda }\,$ and 
$\,\Phi_{s,\mu}\,$ in $\, G \,$ are the same as the
relations satisfied by the corresponding automorphisms
of the {\it free} associative algebra $\, 
\c\langle x,y \rangle\,$. In \cite{BW} this problem
did not arise, because we knew in advance that the map
$\, \omega:\, \mathfrak{C} \to \mathfrak{R} \,$
was bijective, so we had only to transfer to 
$\,  \mathfrak{C} \,$ the natural action of $\, G \,$
on $\, \mathfrak{R} \,$.
\subsection*{}

As explained in the Introduction, to prove Theorem~\ref{IT2}
we have now only to check that the maps $\, \theta \,$
and $\, \omega^{-1} \,$ agree for one point in each $ G$-orbit 
$\, \omega(\mathfrak{C}_n) \subset \mathfrak{R} \,$.
A suitable point is (the class of) the ideal 
$\, I = y^{n+1}A \,+\, (yx-n)A \,$; this is the formal 
Fourier transform of the ideal $\, M \,$ in Section~\ref{EX}. 
As in \cite{BW}, 
we identify $\, A \,$ with the ring 
$\,\c[z, \partial/\partial z] \,$ of differential
operators with polynomial coefficients by
$\, x \leftrightarrow \partial/\partial z \,$,
$\, y \leftrightarrow  z \,$. Then 
$\, I \cap \c[z] \not= \{0\} \,$, so we can
calculate that the Cannings-Holland map 
sends $\, I \,$ to the point 
$$
W = z^{-1} \left\{\, f \in \c[z]\ |\  f^{(n)}(0) = 0 \,\right\}
\ \in \ \mbox{\rm Gr}^{\mbox{\scriptsize \rm ad}}\ .
$$
The (reduced stationary) Baker function of this point is
$$
\tilde{\psi}_{W}(x,z) = 1 - n \, x^{-1} z^{-1}\ .
$$
If $\, \X_n \,$, $\, \Y_n \,$ are the two matrices
found in the example of Section~\ref{EX}, then we have
$$
\tilde{\psi}_{W}(x,z) = \det\left\{\,\I - 
(x\I - \Y_n)^{-1} (z\I + \X_n)^{-1}\,\right\}\ ,
$$
which means that the map $\, \omega^{-1} \,$ sends $\, I \,$
to the pair of matrices $\,(\Y_n, \, -\X_n)\,$. 
The Fourier transform on $ \mathfrak{R} $ corresponds
to the map $\, (\X,\,\Y) \mapsto (-\Y, \, \X)\,$ on
matrices, hence indeed $\, \omega^{-1} \,$ sends the
ideal $\, M \,$ to $\,(\X_n, \, \Y_n) = \theta(M)\,$.

\section{The Beilinson Equivalence}
\la{Sect8}

For each $\, i \in \{0,1, \ldots,|\bs{w}|\} \,$,  
we set $\, \ms{E}_{i} := \ms{O}_{\P(\bs{w})}(i) \,$, 
and $\, \ms{E} := \bigoplus_{i=0}^{|\bs{w}|} 
\ms{E}_{i}\,$. Let
$$
B := \mbox{\rm Hom}(\ms{E},\ms{E}) = \bigoplus_{i,j=0}^{|\bs{w}|} 
\mbox{\rm Hom}(\ms{E}_{i},\ms{E}_{j}) 
$$
be the algebra of endomorphisms of $\, \ms{E}\,$. We consider the 
(left exact) functor $\, \mbox{\rm Hom}(\ms{E},\mbox{---})\,$,
which takes (quasi)coherent sheaves over $\, \P(\bs{w}) \,$ to
right $ B$-modules. Using the fact that $\, \P(\bs{w}) \,$
has finite cohomological dimension and the 
Finiteness Theorem~\ref{FVT}$(a)$, we see that 
$\, \mbox{\rm Hom}(\ms{E},\,\mbox{---}\,)\,$ extends to 
a functor on bounded derived categories 
\begin{equation}
\la{6.16}
\mbox{\rm \textbf{R}Hom}(\ms{E},\, \mbox{---}\,)\,:\,
D^{b}({\sf coh}\,\P) \to 
D^{b}({\sf mod}\,B)\ ,
\end{equation}
where $\, {\sf mod}(B) \,$ denotes the category of 
finite-dimensional right $ B$-modules. 
The following statement is the analogue for $\,\P(\bs{w})\,$
of a theorem of Beilinson (see \cite{B}) for the usual 
projective spaces $\,\mathbb{P}^n \,$.
\begin{theorem}
\la{Beil}
The functor~(\ref{6.16}) is an equivalence of categories. 
\end{theorem}
\begin{proof}
According to \cite{Bon}, Theorem~6.2, it is enough to check that
the sequence of sheaves $\, (\ms{E}_{0}, \ms{E}_{1}, 
\ldots, \ms{E}_{|\bs{w}|}) \,$ (regarded as $0$-complexes 
in $\,D^{b}({\sf coh}\,\P)$) is a {\it complete strongly 
exceptional collection}. Here ``complete'' means that these
objects generate $\,D^{b}({\sf coh}\,\P)\,$ as a triangulated 
category, while ``strongly exceptional'' means that
$$
\begin{array}{lll}
(a) & \mbox{Hom}(\ms{E}_{i}, \ms{E}_{j}) = 0 & 
\quad \mbox{if}\ i > j \ ,\\*[1ex]
(b) & \mbox{\rm Ext}^{k}(\ms{E}_{i},\ms{E}_{j}) = 0 
& \quad \mbox{for all}\ i,j\in \{ 0,1, \ldots, |\bs{w}|\} 
\ \mbox{and}\ k \not= 0\ .
\end{array}
$$
Property $(a)$ is trivial, since
$$
\mbox{Hom}(\ms{E}_{i},\ms{E}_{j}) \cong
\mbox{Hom}(\ms{O}_{\P}(i), \ms{O}_{\P}(j)) \cong
H^{0}(\P(\bs{w}),\, {\mathcal O}_{\P}(j-i)) 
\cong A_{j-i}(\bs{w})\ . 
$$
Similarly, $\, \mbox{\rm Ext}^{k}(\ms{E}_{i},\ms{E}_{j}) \cong
H^{k}(\P(\bs{w}),\, {\mathcal O}_{\P}(j-i)) \,$. 
By Theorem~\ref{ASL} and Proposition~\ref{L4.1} we have 
$\, H^1(\P,\, \ms{O}_{\P}(r)) = 0\, $ and 
$\, H^2(\P,\, \ms{O}_{\P}(r)) \cong (A_{r-|\bs{w}|-1})^{*} \,$ 
for all $\, r \,$. If $\,i,j\in \{ 0,1, \ldots, |\bs{w}|\}\,$ 
then $\, j-i \leq |\bs{w}|\,$, hence 
$\, (j-i)-|\bs{w}|-1 < 0 \,$. Property $(b)$ above follows.

It remains to show that the collection $\, (\ms{E}_{0}, \ms{E}_{1}, 
\ldots, \ms{E}_{|\bs{w}|}) \,$ is complete. 
Denote by $\, E \,$ the smallest
strictly full triangulated subcategory of $\,D^{b}({\sf coh}\,\P)\,$
containing the objects $\, \ms{O}_{\P}, \ms{O}_{\P}(1), \ldots, 
\ms{O}_{\P}(|{\boldsymbol w}|) \,$.  We must show that 
$\, E = D^{b}({\sf coh}\,\P)\, $.
Since any derived category is generated by its abelian 
core, it suffices to prove that the $0$-complexes
$\, (\,\cdots  \to  0 \to \ms{M} \to 0 \to \cdots\,) \, $ are in $\,E\,$ for all
$\, \ms{M} \in \mbox{\sf coh}(\P) \,$.
We know that $\, \bs{A} \,$ has finite global dimension, 
so every $\, \bs{M} \in {\sf grmod}(\bs{A}) \,$ has
a finite projective resolution. Moreover, 
every graded projective $\bs{A}$-module is a finite direct sum 
of shifts of $\, \bs{A} \,$ (see \cite{CE}, Theorem~6.1).
Therefore every $\, \ms{M} \in \mbox{\sf coh}(\P) \,$
has a finite resolution by finite direct sums of
sheaves $\, \ms{O}_{\P}(m)\, $. Such a resolution
gives a complex isomorphic to $ \ms{M} $ in the derived category, 
and hence  (see \cite{GM}, Chapter~III, 
\S~5, Exercise~4(b)), we need only to show that 
$\, \ms{O}_{\P}(m)\,$ belongs to  $\, E \,$ 
for any $\, m \in \Z\,$. 

According to \cite{Steph} (combine Proposition~2.5(ii) and 
Corollary~2.6(ii)), the trivial $ \bs{A}$-module 
$\, \bs{A}/\bs{A}_{\geq 1} \cong \c \,$ has a graded 
resolution of the form
\begin{eqnarray}
\la{6.6}
\lefteqn{
0 \to \bs{A}(-|{\boldsymbol w}|-1) \to 
\bs{A}(- w_2 - 1) \oplus \bs{A}(- w_1 - 1) 
\oplus \bs{A}(-|{\boldsymbol w}|) \to} \\*[2ex] 
& & \to \bs{A}(- w_1) \oplus \bs{A}(- w_2) 
\oplus \bs{A}(-1) \to 
\bs{A} \to {\mathbb C} \to 0\ . \nonumber
\end{eqnarray}
By shifting degrees in (\ref{6.6}) and passing to the quotient
category we get (for any integer $\, m $) 
a Koszul-type exact sequence in 
$\, \mbox{\sf coh}({\mathbb P}_{q}^{2}) \,$:
\begin{eqnarray}
0 & \to & {\mathcal O}(m) \to 
          {\mathcal O}(m+w_1) \oplus {\mathcal O}(m+w_2) \oplus 
          {\mathcal O}(m+1) \to \nonumber \\*[2ex]
  & \to & {\mathcal O}(m+w_2+1) \oplus {\mathcal O}(m+w_1+1) 
          \oplus {\mathcal O}(m+|{\boldsymbol w}|) \to 
          {\mathcal O}(m+|{\boldsymbol w}|+1) \to 0\ . \nonumber
\end{eqnarray}
Letting $\, m = 0 \,$ above, we observe that 
$\,\ms{O}(|{\boldsymbol w}|+1) \,$ is quasi-isomorphic to a complex
each term of which is in $\, E \,$; therefore, 
$ \, {\mathcal O}(|{\boldsymbol w}|+1) \in E \,$.
Similarly, for $ m = -1 $, it follows that
$ \, {\mathcal O}(-1) \in  E\,$. Arguing in 
this way by induction (going both in negative and positive directions), 
we conclude that $ \, \ms{O}(m) \in E \,$ 
for all $\, m \in {\mathbb Z}\,$. This finishes the proof 
of the theorem.
\end{proof}
We can now explain the significance of the vanishing conditions 
(\ref{I2}) (and their generalization in Theorem~\ref{PP1}).
Given a coherent sheaf $\, \ms{M}\,$ over $\, \P \,$, we may regard it as
a $0$-complex in $\, D^{b}({\sf coh}\,\P) \,$. By definition,
the functor $\, \mbox{\rm \textbf{R}Hom}(\ms{E},\, \mbox{---}\,) \,$
then maps $\, \ms{M} \,$ to a complex of $ B$-modules whose
cohomology in degree $\, -k \,$ is
$$
\mbox{\rm Ext}^{k}(\ms{E},\, \ms{M})\,
\cong \,\bigoplus_{i=0}^{|\boldsymbol{w}|} 
\mbox{\rm Ext}^{k}(\ms{O}_{\P}(i),\, \mathcal{M}) \,\cong\,
\bigoplus_{i=0}^{|\boldsymbol{w}|}
H^{k}(\P,\, \mathcal{M}(-i))\ .
$$
In our case Theorem~\ref{PP1} tells us that 
this cohomology vanishes for all $ \, k \not= 1\,$. 
A complex with cohomology
only in one degree is isomorphic (in the derived category)
to its cohomology; thus, essentially, the Beilinson equivalence
assigns to our sheaf $\, \ms{M} \,$ the single $ B$-module 
$\, \bigoplus_{i=0}^{|\boldsymbol{w}|} 
H^{1}(\P,\, \mathcal{M}(-i))\,$.

\begin{figure}
\centering
\vspace{1.3in}
\hspace{-1.6in}
\includegraphics[totalheight=2in]{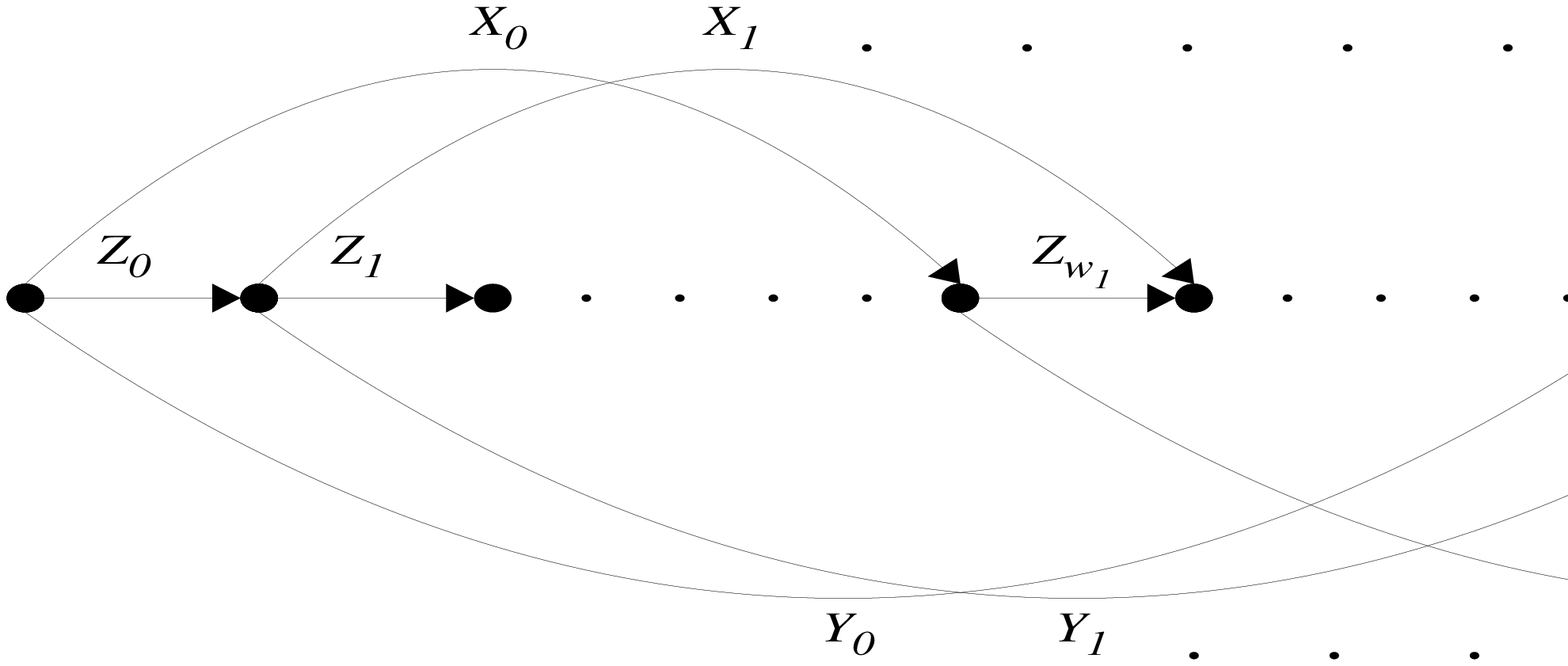}
\vspace{-1.2in}
\caption{}
\label{Fig1}
\end{figure}

To make contact with the language of quivers used 
in the Introduction, we have only to note that 
the algebra $\, B \,$ is isomorphic to the path algebra
of the quiver\footnote{In the case $\, \bs{w} = (1,1) \,$, this is
just the quiver indicated in (\ref{I3}).} 
shown in Fig.~\ref{Fig1},
%
%
%
with the relations
\begin{equation}
\la{qrel}
\begin{array}{lcl}
Z_{i+ w_1}\, X_{i} - X_{i+1}\, Z_{i} & = & 0 
\ , \qquad i=0,1, \ldots, w_2-1\ ; \\*[1ex]
Z_{j+ w_2}\, Y_{j} - Y_{j+1}\, Z_{j}  & = & 0 \ , 
\qquad j =0,1, \ldots, w_1-1\ ; \\*[1ex]
X_{w_{2}} \, Y_{0} \, - \, Y_{w_{1}}\, X_{0} & = & 
Z_{|{\boldsymbol w}|-1}\, Z_{|{\boldsymbol w}|-2}\, 
\ldots \, Z_{1}\, Z_{0}\ .
\end{array}
\end{equation}
That means that a (left) $ B$-module can be identified
with a representation of this quiver with relations.
Indeed, let $\, e_{i} \,$ denote the identity map in
$\, \mbox{Hom}(\ms{E}_i, \ms{E}_i) \subset B \,$:
these are mutually orthogonal idempotents in 
$\, B \,$, and $\, e_0 + e_1 + \ldots + 
e_{|\bs{w}|} = 1_{B} \,$. Hence any $ B$-module
$\, V \,$ decomposes: $\, V = \bigoplus_{i} V_i \,$,
where $\, V_i := e_i V \,$. Further, each
element of $\, \mbox{Hom}(\ms{E}_i, \ms{E}_j) 
\subset B \,$  maps $\, V_i \,$ to $\, V_j \,$. 
As we have seen above, 
$\, \mbox{Hom}(\ms{E}_i, \ms{E}_j) \,$ is 
naturally isomorphic to  
$\, A_{j-i}(\bs{w})\, $; hence the generators
$\, X, \, Y,\,$ and $\, Z \,$ of $\, \bs{A} \,$
determine maps $\, X_{\alpha}\, , \, Y_{\alpha}\,$ and 
$\, Z_{\alpha} \,$ as indicated in Fig.~\ref{Fig1}
(the vertices in this diagram represent the spaces
$\, V_i $). The relations (\ref{qrel}) follow from the
defining relations (\ref{4.4}) of the algebra $\, \bs{A}(\bs{w}) \,$.
In this way, each $B$-module determines a representation of 
the above quiver with relations. The construction of a
 $B$-module from a quiver representation is equally 
straightforward.

\section{Le Bruyn's Moduli Spaces}
\la{Sect9}
As we mentioned in the Introduction, our construction
differs from that in Le Bruyn's paper \cite{LeB} only
in the different choice of extension of an $ A$-module
$ M $ to a sheaf over $ \P \,$. In this section we 
clarify the relationship between the two constructions.
We confine ourselves to the (basic) case $\, \bs{w} = (1,1)\,$.

As usual, let $ M $ be an ideal of $ A\,$, with the 
normalized induced filtration, and 
let $\, d \,$ be the minimum filtration degree of elements 
of $ M\,$. We recall (see Lemma~\ref{min1}) that $\, d \geq 1\,$
(unless $ M $ is cyclic: this case has to be 
excluded from some of the statements below).
It follows from Theorem~\ref{PP1} that the sheaf 
$\, \ms{F} = \ms{M}(d-1) \, $
satisfies the vanishing conditions (\ref{I2}), so that
$\, \ms{M}(d-1) \, $ (and hence $ M $) is determined by 
the quiver representation
\begin{equation}
\label{9.1}
H^{1}(\P, \ms{M}(d-3))\  
\begin{array}{c}
\longrightarrow \\*[-1ex] 
\longrightarrow \\*[-1ex]
\longrightarrow 
\end{array} 
\ H^{1}(\P,  \ms{M}(d-2))\  
\begin{array}{c}
\longrightarrow \\*[-1ex] 
\longrightarrow \\*[-1ex]
\longrightarrow 
\end{array}
\ H^{1}(\P,  \ms{M}(d-1)) \ .
\end{equation}
Le Bruyn uses a slightly more subtle fact 
(see \cite{Baer}, Corollary~7.2): because  
the sheaf $\, \ms{F} = \ms{M}(d-2) \, $ also 
satisfies (\ref{I2}), $ M $ is determined
by the left hand part 
$$
H^{1}(\P,  \ms{M}(d-3))\  
\begin{array}{c}
\longrightarrow \\*[-1ex] 
\longrightarrow \\*[-1ex]
\longrightarrow 
\end{array}
\ H^{1}(\P, \ms{M}(d-2)) 
$$
of (\ref{9.1}). For each pair of non-negative integers
$ (r,s) \,$, let $\, \tilde{\mathfrak{M}}(r,s) \,$ be the space
of isomorphism classes of quintuples $\, 
(V, \, W; \, X_1, \, X_2, \,X_3) \,$, 
where $ V $ and $ W $ are vector spaces of 
dimensions $ r $ and $ s $ (respectively) and
$\, X_i \,$ are linear maps from  $ V $ to $ W \,$.
Thus each ideal  $ M $ of $ A $ determines (and is determined by)
a point of $\, \tilde{\mathfrak{M}}(r,s) \,$, where
$\, r = \dim_{\mathbb{C}} H^{1}(\P,  \ms{M}(d-3))\,$ and
$\, s = \dim_{\mathbb{C}} H^{1}(\P,  \ms{M}(d-2))\,$.
Denoting by $\, \mathfrak{M}(r,s) \,$  the subset of points
of $\, \tilde{\mathfrak{M}}(r,s) \,$ that arise in this way,
we obtain the main result of \cite{LeB}: the space of ideals
$\, \mathfrak{R}\, $ decomposes as the disjoint union of the 
``moduli spaces'' $\, \mathfrak{M}(r,s) \,$ ($ r,s \geq 0 $).

The relationship of this to our decomposition of 
$\, \mathfrak{R}\, $ becomes clear if we calculate
the dimensions $ r $ and $ s $ in terms of $ d $ 
and our invariant $ n \,$. Because of 
(\ref{I2}), $\, -r \,$ and $\, -s \,$  are equal
to the Euler characteristics $\, \chi(\P,\,\ms{M}(d-3))\,$ and 
$\, \chi(\P,\, \ms{M}(d-2)) \,$ respectively.
From (\ref{DIM2}) we get
$$
\chi(\P,\,\ms{M}(k)) = \frac{1}{2}(k+1)(k+2) - n \quad \mbox{for} \ 
k \geq -2 \ 
$$ 
(in fact, this formula is true for all $\, k \in \Z$).
In particular, we have
\begin{equation}
\la{9.3}
r = n - \frac{1}{2}\,(d-2)(d-1)\quad , \quad
s = n - \frac{1}{2}\,d(d-1)\ .
\end{equation}
From this we notice immediately that the space 
$\, \mathfrak{M}(r,s) \,$ is empty if $ r < s\,$,
while in general each point of $\, \mathfrak{M}(r,s) \,$
determines an ideal $ M $ for which our parameters are 
given by
\begin{equation}
\label{9.4}
n = \frac{1}{2}\left[(r-s)^2 + r + s\right]
\quad , \quad  d = r - s + 1\ .
\end{equation}
The map $\, \theta \,$ then gives us a point of 
$\, \mathfrak{C}_{n} \,$. We thus have
\begin{proposition}
\la{P9.1}
Let $\,\mathfrak{C}_{n}(d)\,$ be the subspace of 
$\,\mathfrak{C}_{n}\,$ corresponding to ideals of 
minimum filtration degree $ d\,$. Then the 
construction explained above defines a bijection
$\, \mathfrak{M}(r,s) \to \mathfrak{C}_{n}(d) \,$,
where the numbers $ (r,s) $  and $ (n,d) $ are 
related by (\ref{9.3}) and (\ref{9.4}).
\end{proposition}

The commutative analogue of this decomposition
of $\,\mathfrak{C}_{n}\,$ 
(into the subspaces $\,\mathfrak{C}_{n}(d)\,$)
forms part of the much studied\footnote{We thank A.~Iarrobino for
information on this subject.} 
{\it Brill-Noether theory}\,: in that case 
$\,\mathfrak{C}_{n}\,$ is replaced by the point
Hilbert scheme $\, \mbox{Hilb}_{n}(\mathbb{A}^2)\,$
and $\,\mathfrak{C}_{n}(d)\,$ by the 
subvariety $\, \mbox{Hilb}_{n}(d)\,$ (say) of 
$n$-tuples of points of $\, \mathbb{A}^2 \,$
that lie on a curve of degree $\, d \,$, but not on one of 
degree $ d-1\,$. The detailed structure of this stratification 
of $\, \mbox{Hilb}_{n}(\mathbb{A}^2)\,$ seems quite 
complicated (see, for example, \cite{BH}, \cite{R}). 
However,  using the fact that there is a curve of degree 
$ d $ through any $\, d(d+3)/2 \,$ points in the plane, it is 
easy to see that (for $ n>0 $)  $\, \mbox{Hilb}_{n}(d)\,$
is non-empty if and only if we have 
$\,1 \leq d \leq D \,$, where $ D $ is the least integer
such that 
$\, n \leq D(D+3)/2 \,$. Furthermore, the dimension
of $\,\mbox{Hilb}_{n}(d)\,$ is then given by 
$$ 
\dim_{\c} \mbox{Hilb}_{n}(D) = 2n \ ,\quad \dim_{\c}
\mbox{Hilb}_{n}(d) = n + \frac{1}{2}\,d(d+3)\quad \mbox{if} \quad
d < D\ .
$$
We expect that the situation is the same in the noncommutative 
case. It might be interesting to study this decomposition
of $\,\mathfrak{C}_{n}\,$ in more detail to see to what
extent it is simply a deformation of what we have in the 
commutative case.

\section{Ideals and Bundles}
\la{Sect10}
In the case $\, \bs{w} = (1,1) \,$, the authors of \cite{KKO} 
establish a bijection between $\, \mathfrak{C} \,$ and the space 
$\, \mathfrak{L} \,$ of all line bundles (suitably defined) 
over $ \P $ that are trivial on the line at infinity. This
bijection is constructed using monads, following the original 
approach of Barth to classifying bundles over projective spaces
(see \cite{Bar}, \cite{N}). Here we shall check that the result of
\cite{KKO} is essentially equivalent to the
bijectivity of our map $\,\theta: \mathfrak{R} \to  \mathfrak{C} \,$.
Most of what follows is valid for any positive weight $\, \bs{w} \,$.
The key step is the following lemma, which may be of independent 
interest.
\begin{lemma}
\la{ind}
Let $ M $ be a finitely generated rank one torsion-free $ A$-module,
and let $\,\{ M_{\bullet}\}\,$ be any filtration of $ M\,$ by 
finite-dimensional subspaces $ M_k\,$. Suppose 
the associated graded module $\, \bs{GM} \,$ is 
{\rm essentially torsion-free} (meaning that $\, \bs{GM}_{\geq N} \,$ 
is torsion-free for some $ N \,$).
Then for any embedding of $ M $ in $ A $ there is an integer
$ k_{0} $ such that $\, M_k = M \cap A_{k-k_0} \,$ for all 
$\, k \geq N \,$.  In other words, the given filtration on $ M $
essentially coincides with an induced filtration.
\end{lemma}
\begin{proof}
Let $ \u $ denote the valuation on $ M $ corresponding
to the given filtration, that is, if $\, m \in M\,$ then
$ \u(m) $ is the smallest integer $ k $ such that $\, m \in M_k\,$.
Let $ \vv $ be the valuation on $ A $ corresponding to 
the ring filtration we are using. Then
\begin{equation}
\la{v1}
\u(ma) \leq \u(m)+ \vv(a) \quad \mbox{for all}\ m \in M\  
\mbox{and}\ a \in A \ .
\end{equation}
The assumption on $\, \bs{GM} \,$ is equivalent to
\begin{equation}
\la{v2}
\u(ma) = \u(m)+ \vv(a) \quad \mbox{if} \quad  \u(m) \geq N\ .
\end{equation}
Fix an embedding of $ M $ in $ A \,$. 
Then $\, M \otimes_{A} Q = Q \,$, hence every element
$ q $ of the Weyl quotient field $ Q $ can be written in the 
form $\, q = m a^{-1}\,$ with $\, m \in M \,$ and $\, a \in A\,$.
Moreover, we can choose $ m $ so that $\, \u(m) \geq N \,$. 
Indeed, for any $\, b \in A \,$ we have 
$\, m a^{-1} = mb (ab)^{-1} \,$, and we cannot have $\, \u(mb) < N \,$ 
for all $ b $ because $ M_N $ is finite-dimensional.

Now define a function $\, F : Q \to \Z \,$ as follows:
if $\, q = m a^{-1}\,$ as above with $\, \u(m) \geq N \,$,
let $ F(q) := \u(m) - \vv(a) \,$. To show that $ F $ is well
defined, suppose that $ m a^{-1} = n b^{-1} $ are two such 
expressions for $ q $. Since $ A $ is an Ore domain, we have
$\, b^{-1} a = p r^{-1}\,$ for some $\, p,r \in A\,$.
So $\, mr = np \,$ and $\,ar = bp\,$. Using (\ref{v2})
and the similar fact for the filtration on $ A \,$,
we get
$$
\u(m)+ \vv(r) = \u(n) + \vv(p) \quad \mbox{and} 
\quad  
\vv(a) + \vv(r) = \vv(b) + \vv(p)\ .
$$
Hence $\, \u(m) - \vv(a) = \u(n) - \vv(b)\,$, as desired. 

Now, if $\, \u(m) \geq N \,$ then 
$\, F(1) = F(mm^{-1}) = \u(m) - \vv(m) \,$. Thus, setting 
$\, k_0 :=  F(1) \,$, we have
\begin{equation}
\la{v3}
\vv(m) = \u(m) - k_0 \quad \mbox{if} \quad  \u(m) \geq  N\ .
\end{equation}
Also, for {\it any} element $\, m \in M\,$, if we choose $\, a \in A\,$
so that $\, \u(ma) \geq N \,$, we have (using (\ref{v1}))
$\,
\vv(ma) = \vv(m)+ \vv(a) = \u(ma) - k_0 \leq \u(m) + \vv(a)- k_0 \ ,
\,$
so
\begin{equation}
\la{v4}
\vv(m) \leq \u(m) - k_0 \quad \mbox{for all} \quad  m \in M\ .
\end{equation}
Let $\, \{ M_{\bullet}'\} \,$ denote the induced filtration on 
$ M \,$. By (\ref{v4}),  we have 
$$
m \in M_{k} \ \Leftrightarrow \ \u(m) \leq k \ \Rightarrow \
\vv(m) \leq k-k_0  \ \Leftrightarrow \ m \in M_{k-k_0}'\ ,
$$
that is, 
\begin{equation}
\la{v5}
M_{k} \subseteq M_{k-k_0}' \quad \mbox{for all}\ k\ .
\end{equation}
Similarly, by (\ref{v3}), if $\, k \geq N \,$, we have
\begin{eqnarray}
\lefteqn{
m \not\in M_{k} \  \Leftrightarrow \  
\u(m) > k \ \Rightarrow \ \u(m) \geq N\ \Rightarrow } \nonumber\\*[0.6ex]
&& \Rightarrow\ \vv(m) = \u(m) - k_0 \ \Rightarrow \ 
\vv(m) > k - k_0 \ \Leftrightarrow\ 
 m \not\in M_{k-k_0}'\ ; \nonumber
\end{eqnarray}
equivalently, $\, m \in M_{k-k_0}' \ \Rightarrow \ m \in M_{k} \,$,
and hence 
\begin{equation}
\la{v6}
M_{k-k_0}' \subseteq  M_{k} \quad \mbox{for all}\ k \geq N\ .
\end{equation}
The lemma now follows at once from (\ref{v5}) and (\ref{v6}).
\end{proof}
\begin{lemma}
\la{bundle}
Let $\, M \,$ be an ideal of $ A \,$, and let $ \, \ms{M}\,$
be its canonical extension to $ \P \,$. Then  
$\, \ms{M}\,$ is a bundle in the sense of \cite{KKO} 
(Definition~5.4).
\end{lemma}
\begin{proof}
According to \cite{KKO} (see Section~5.3), we have to show 
that 
$$
\mbox{\rm Ext}^{i}(\ms{M}, \ms{O}(k)) = 0 \quad \mbox{for}\  
k \gg 0 \  \mbox{and} \ i > 0\ .
$$
By Serre Duality (see Theorem~\ref{FVT}$(c)$), it is equivalent
to show that $\, H^{0}(\P\,,\,\ms{M}(k)) \,$ and
$\, H^{1}(\P\,,\,\ms{M}(k)) \,$ vanish for $\, k \ll 0\,$. 
The vanishing of $\,H^{0}\,$ is part of Theorem~\ref{PP1}$(ii)$.
The statement about $\,H^{1}\,$ follows from 
Theorem~\ref{PP1} $(i)$ by a useful argument which would not
be available in the commutative case. 
Because the spaces $\, H^{1}(\P\,,\,\ms{M}(k)) \,$ 
are finite-dimensional,  Theorem~\ref{PP1} $(i)$ implies 
that the maps $\, H^{1}(\P\,,\,\ms{M}(k-1)) 
\to H^{1}(\P\,,\,\ms{M}(k)) \,$  induced by 
multiplication by $\, Z \,$ are 
isomorphisms for $\, k \ll 0\,$. 
If we use these isomorphisms to identify the spaces 
$\, H^{1}(\P\,,\,\ms{M}(k)) \,$ for $\, k \ll 0\,$,
the action of the generators $ X $ and $ Y $ gives us
a finite-dimensional representation of the Weyl algebra,
which is impossible unless the representation space is zero. 
\end{proof}
\begin{proposition}
\la{monad}
Let $\, \mu: \mathfrak{R} \to  \mathfrak{L}  \,$
be the map that sends an ideal class to its
canonical extension. Then $\, \mu \,$ is bijective.
\end{proposition}
\begin{proof}
The inverse map $\, \nu: \mathfrak{L} \to \mathfrak{R}\,$ 
is constructed as follows.  
Let $\, \ms{M} \,$ be a line bundle 
over $ \P \,$, trivial over $\, l_{\infty}\,$, and let $\, 
\bs{M} = \oplus\,M_k \,$ be 
a graded $ \bs{A}$-module with $ \pi\bs{M} = \ms{M} \,$.
By (the proof of) Lemma~6.1 of \cite{KKO}, $\, \ms{M} \,$
can be embedded in a direct sum of sheaves 
$\, \ms{O}(k) \,$; hence the $ \bs{A}$-module
$\, \bs{M} \,$ is essentially torsion-free.
We set $\, M := \lim\limits_{\longrightarrow}^{} M_k \,$,
where the direct limit is taken over the maps $\, \cdot Z\,: \,
M_{k-1} \to M_k \,$. Then $\, M \,$ is a rank one torsion-free 
$A$-module, filtered by (the images of) the components $\, M_k \,$.
Forgetting this filtration, we obtain a map 
$\, \nu: \mathfrak{L} \to \mathfrak{R}\,$.
Since $\, \ms{M} \,$ is trivial over $\, l_{\infty} \,$, we have 
$\, \bs{GM}_{\geq k} \cong \bs{GA}_{\geq k}\,$ for $\, k \gg 0\,$, 
and therefore the filtration on $ M $ coincides 
(in sufficiently high degrees) with the normalized induced 
filtration (see Lemma~\ref{ind}). It follows easily that 
$\,\nu = \mu^{-1}\,$.
\end{proof}

We omit the proof that the bijection $\,\theta\,\nu: \mathfrak{L} \to
\mathfrak{C} \,$ coincides\footnote{Actually, there is a difference
of sign.} with the map constructed in \cite{KKO}. Although we do not
know a reference for this fact, it is very unsurprising, since
Beilinson's equivalence is in essence a generalization of the 
monad construction (cf. \cite{B}).

\appendix
\section{Appendix by Michel Van den Bergh}
In this appendix we give alternative proofs of Theorems 1.3 and 1.4.
Our proof of Theorem 1.4
does not rely on $\check{\mathrm{C}}$ech cohomology. Furthermore our
proof of Theorem 1.3 does not rely on the properties  of weighted projective
spaces with $w\neq (1,1)$. So it is in fact independent of Theorem
1.4!

After the authors of this paper had proposed me to write this appendix 
they succeeded in simplifying some of my original arguments and they
have gracefully allowed me to consult some of their private notes
which contained similar ideas. These contributions have allowed me to
streamline the presentation below. 

The main idea behind the new proofs is that while the map (a priori
dependent on $w$) which
associates linear data to ideals  seems hard to understand, the inverse
of that map is given by a simple formula (see \eqref{leftinverse1}
below) whose properties are transparent.

We start with the proof of Theorem 1.4. First we introduce some notational
conventions. If $Q$ is a quiver with relations then we will
identify $Q$ with a $\CC$-linear additive category whose objects are
the finite
direct sums of vertices of $Q$. This has the effect that we write paths from
right to left. Under this formalism the path algebra $\CC Q$ of $Q$ is given
by the endomorphism ring of the sum of the vertices. Note that $\CC
(Q^{\opp})=(\CC Q)^{\opp}$. A morphism of quivers is  a functor between the associated
categories. Such a morphism induces a ring homomorphism between the
associated path algebras. 

If $Q$ is a quiver with relations then $\Mod(Q)$ is the category of $\CC$-linear contravariant functors on $Q$
with values in $\CC$-vector spaces.   By Yoneda's lemma we obtain a
full faithful functor $Q\r \Mod(Q)$ whose image consists of finitely
generated projectives. 
Invoking Morita theory or
directly one sees that $\Mod(Q)$ is equivalent to $\Mod(kQ)$, the
category of right $kQ$-modules.

We now use the notations from Section  8. We denote the quiver (with relations) given in Figure 1 by
$\Delta$. If we view $\Delta$ as a $\CC$-linear additive
category then it is equivalent to the
full subcategory of $\coh \PP^{2}_{q}$ whose objects are finite direct
sums of the $\Oscr_{\PP^2_q}(i)_{i=0,\ldots
  |w|}$ in such a way that the $i$'th
vertex from the left (counting from $0$) corresponds to
$\Oscr_{\PP^2_q}(i)$. It follows that $\CC
\Delta=\End(\oplus_{i=0}^{|w|}\Oscr_{\PP^2_q}(i))=\End(\Escr)=B$.

As noted in Section  8, the functor $\RHom_{\PP^2_w}(\Escr,-)$
defines an equivalence between the triangulated categories $D^b_f(\coh(\PP^2_q))$
and $D^b(\mod(B))$. The inverse
functor is given by $-\Lotimes_{B} \Escr$. It is clear that this
equivalence restricts to an equivalence between the following two
subcategories 
\[
\Xscr_{1}=\{\Mscr\in \coh(\PP^2_q)\mid \Ext^i_{\PP^2_q}(\Escr,\Mscr)=0\text{ for }i\neq 1
\}
\]
and
\[
\Yscr_{1}= \{M\in \mod(B)\mid \Tor_i^{B}(M,\Escr)=0\text{ for }i\neq 1
\}
\]
The inverse equivalences between these categories are given by
$\Ext^1_{\PP^2_q}(\Escr,-)$ and $\Tor^{B}_1(-,\Escr)$.

As before we denote by $\PP^1$ the line at infinity in $\PP^2_q$. Note
that  $\PP^1$ is a weighted projective line in
the sense of \cite{GL}. The inclusion $\PP^1\r \PP^2_q$ is denoted by
$i$. 

Let us denote by $\Rscr$ the full subcategory of $\coh(\PP^2_q)$ whose
objects have the property
that $\Mscr\not\cong \Oscr_{\PP^2}$ and $i^\ast(\Mscr)\cong
\Oscr_{\PP^1}$.
Using the results in Section  4 one shows that $\Rscr\subset \Xscr_1$ and
furthermore that the image of $\Rscr$ under $\Ext^1(\Escr,-)$ lies in
the following category
\begin{multline*}
\Cscr_{1}=\{M \in \Mod(\Delta)\mid M(Z_i)\text{ is an isomorphism
  for }i=1,\ldots,|w|-1,\\\text{ $M(Z_{0})$ is surjective with one
  dimensional kernel and all $M(i)$}\\
\text{ are finite dimensional}
\}
\end{multline*}
We define $\Cscr_2$ as the category consisting of triples $(W,\XX,\YY)$
where $W$ is a finite dimensional vector space and $\XX,\YY$ are
endomorphisms of $W$ satisfying $\rk ([\YY,\XX]-\Id)=1$. 
Let $M\in \Cscr_{1}$. Then up to a canonical isomorphism we may
assume that $M(1)=\cdots=M(|w|)$ and
$M(Z_1)=\cdots=M(Z_{|w|-1})=\Id$. 

Put $W=M(|w|)$, $\XX=M(X_1)$,
$\YY=M(Y_1)$ and 
$\ZZ=M(Z_{0})$. Taking into account that $M$ is a contravariant
functor we find 
\begin{gather*}
M(X_1)=M(X_2)=\cdots=M(X_{w_2})=\XX\\
M(Y_1)=M(Y_2)=\cdots=M(Y_{w_1})=\YY\\
M(X_0)=\ZZ\XX\\
M(Y_0)=\ZZ\YY
\end{gather*}
It follows that $\ZZ(\YY\XX-\XX\YY-\Id)=0$ and hence $(W,\XX,\YY)\in \Cscr_2$. It is
clear that this procedure is reversible and defines an equivalence
$\Cscr_2\cong \Cscr_{1}$.

Let $R$ be the category of non-trivial rank one projective right
$A$-modules (with maps given by isomorphisms). If 
$N\in R$ then according to Section  4 there
exists, up to isomorphism, a unique extension $\Mscr$ of $N$ to $\PP^2_q$ which
lies in $\Rscr$.

Summarizing we now have a composition of functors:
\begin{equation}
\label{app1}
R\xrightarrow{\cong} \Rscr \hookrightarrow \Cscr_1\cong \Cscr_2
\end{equation}
(that the first functor is an equivalence follows for example
from the easily proved fact that the objects in $\Cscr_2$ are simple
objects when considered as representations of the two loop quiver).
\begin{lemma} \eqref{app1} is an equivalence.
\end{lemma}
\begin{proof} Note that we are not allowed deduce this result from Theorem 1.1 since the
  proof of that theorem depends on Theorem 1.3! 

Below we will construct a (left) inverse to
  \eqref{app1} which is independent of $w$. This means that in principle we have to prove the
  lemma only for one particular $w$. If $w=(1,1)$ then the lemma
  can  be deduced from the results in
  \cite{Ginzburg}, \cite{KKO}, \cite{LeB} although the point of view in these
  papers is slightly different.

We will give a proof which works equally well for all $w$. Perhaps the
method has some independent interest.

 Let $M\in
\Cscr_{1}$. We need to  prove two things:
\begin{enumerate}
\item $M\in \Yscr_1$, i.e. $M\Lotimes_{B} \Escr$ has its only
  non-vanishing cohomology in degree -1. This has the effect that $M$
  is in the image of some object in $\Xscr_1$.
\item $M$ is actually in the image of $\Rscr$, i.e. $i^\ast(H^{-1}(M\Lotimes_{B} \Escr))\cong
  \Oscr_{\PP^1}$. 
\end{enumerate}
Now $\PP^2_q$ has the pleasant property that if $0\neq \Mscr\in
\coh(\PP^2)$ then $i^\ast(\Mscr)\neq 0$. From this we
easily deduce that 1.,2. above are actually equivalent to the following
single statement:
\begin{enumerate}
\setcounter{enumi}{2}
\item $Li^\ast(M\Lotimes_{B} \Escr)\cong \Oscr_{\PP^1}[1]$.
\end{enumerate}
Let $\Escr_\infty=\oplus_{i=0}^{|w|-1} \Oscr_{\PP^1}(i)$ and
$B_{\infty}= \End(\Escr_\infty)$. Then $B_{\infty}$ is the path
algebra of the quiver $\Delta_{\infty}$ 

\vspace*{1.1cm}

\[
\psset{arrows=->,nodesep=0.1cm,labelsep=0.1cm,border=1mm}
\arraycolsep=0.4cm
\begin{array}{ccccccccc}
\rnode{a}{0}  & & \rnode{b}{w_1-1}  & \rnode{c}{w_1}
& &
\rnode{d}{w_2-1} & \rnode{e}{w_2}
& &
\rnode{f}{|w|-1} 
\end{array}
\everypsbox{\scriptstyle}
\ncarc[arcangle=30]{a}{c}
\Aput{X_0}
\ncarc[arcangle=30]{d}{f}
\Aput{X_{w_2-1}}
\ncarc[arcangle=-30]{a}{e}
\Bput{Y_0}
\ncarc[arcangle=-30]{b}{f}
\Bput{Y_{w_1-1}}
\psset{arrows=-, linestyle=dotted}
\ncline{a}{b}
\ncline{c}{d}
\ncline{e}{f}
\]

\vspace*{1.8cm}

Observing that $\RHom(\Escr_\infty,-)$ defines an equivalence between
  $D^b(\coh(\PP^1))$ and $D^b_f(B_\infty)$
we want to understand the composition
\begin{equation}
\label{composition}
D^b_f(B)\xrightarrow{-\Lotimes \Escr} D^b(\coh(\PP^2_q))\xrightarrow{Li^\ast} D^b(\coh(\PP^1))
\xrightarrow{\RHom(\Escr_\infty,-)}D^b_f(B_\infty)
\end{equation}
Checking on projectives, and then on complexes of projectives we find
that on an object $M$ in $\Mod(\Delta)$ the composition
\eqref{composition}  is 
given by a length two complex concentrated in degrees $-1,0$ of the following form

\vspace*{1.1cm}

\[
\psset{arrows=->,nodesep=0.1cm,labelsep=0.1cm,border=0.5mm}
\begin{array}{r@{\hskip 1cm}ccccccccc}
\text{degree 0:}& \rnode{a}{M(0)}  & & \rnode{b}{M(w_1-1)}  & \rnode{c}{M(w_1)}
& &
\rnode{d}{M(w_2-1)} & \rnode{e}{M(w_2)}
& &
\newbox\aa\rnode{f}{\global\setbox\aa\hbox{$M(|w|-1)$}\copy\aa} 
\\[3cm]
\text{degree -1:}&\rnode{a1}{M(1)}  & & \rnode{b1}{M(w_1)}  & \rnode{c1}{M(w_1+1)}
& &
\rnode{d1}{M(w_2)} & \rnode{e1}{M(w_2+1)}
& &
\rnode{xxx}{\hbox to \wd\aa{\hfill \rnode{f1}{$M(|w|)$}\hfill}}
\end{array}
\everypsbox{\scriptscriptstyle}
\ncarc[arcangle=-20]{c}{a}
\Bput{M(X_0)}
\ncarc[arcangle=-20]{f}{d}
\Bput{M(X_{w_2-1})}
\ncarc[arcangle=10]{e}{a}
\aput{0}(0.7){M(Y_0)}
\ncarc[arcangle=10]{f}{b}
\aput{0}(0.35){M(Y_{w_1-1})}
\psset{arrows=-, linestyle=dotted}
\ncline{a}{b}
\ncline{c}{d}
\ncline{e}{f}
\psset{arrows=->, linestyle=solid}
\ncarc[arcangle=-20]{c1}{a1}
\bput{0}(0.65){M(X_1)}
\ncarc[arcangle=-20]{f1}{d1}
\bput{0}(0.35){M(X_{w_2})}
\ncarc[arcangle=10]{e1}{a1}
\aput{0}(0.6){M(Y_1)}
\ncarc[arcangle=10]{f1}{b1}
\aput{0}(0.4){M(Y_{w_1})}
\psset{arrows=-, linestyle=dotted}
\ncline{a1}{b1}
\ncline{c1}{d1} 
\ncline{e1}{f1}
\psset{arrows=->, linestyle=dotted}
\ncline{a1}{a}
\Bput{M(Z_0)}
\ncline{b1}{b}
\Bput{M(Z_{w_1-1})}
\ncline{c1}{c}
\Bput{M(Z_{w_1})}
\ncline{d1}{d}
\Bput{M(Z_{w_2-1})}
\ncline{e1}{e}
\Bput{M(Z_{w_2})}
\ncline{f1}{f}
\Bput{M(Z_{|w|-1})}
\]

\vspace*{1.8cm}

It is now clear that if $M\in \Cscr_{1}$ then the image of $M$ under
the composition \eqref{composition} is equal to
 $S[1]$ where $S$ is the simple object in $\Mod(\Delta_\infty)$ 
defined by $\dim S(i)=\delta_{i0}$. Since $S$ corresponds to $\Oscr_{\PP^1_w}$ we are done.
\end{proof}

Now we continue with the proof of Theorem 1.4.
Theorem 1.4 asserts that the functor \eqref{app1}  is
independent of $w$. We prove this by showing that the inverse of
\eqref{app1} is independent of $w$.

Let $W=(W,\XX,\YY)\in
  \Cscr_2$. Then the associated  object $M$ of $\Mod(\Delta)$  looks
  like

\vspace*{1.1cm}

\begin{equation}
\label{Mdescription}
\psset{arrows=->,nodesep=0.1cm,labelsep=0.1cm,border=1mm}
\arraycolsep=0.4cm
\begin{array}{cccccccccccc}
\rnode{a}{W'}  & \rnode{b}{W} & \rnode{j}{W} & &\rnode{c}{W} & \rnode{d}{W} & &\rnode{e}{W} & \rnode{f}{W} & &
\rnode{g}{W}
& \rnode{h}{W}
\end{array}
\everypsbox{\scriptstyle}
\ncarc[arcangle=-45]{c}{a}
\bput{0}(0.6){\ZZ\XX}
\ncarc[arcangle=-45]{d}{b}
\Bput{\XX}
\ncarc[arcangle=-45]{h}{e}
\Bput{\XX}
\ncarc[arcangle=45]{e}{a}
\aput{0}(0.6){\ZZ\YY}
\ncarc[arcangle=45]{f}{b}
\aput{0}(0.45){\YY}
\ncarc[arcangle=45]{h}{c}
\Aput{\YY}
\ncline{b}{a}
\Bput{\ZZ}
\ncline{j}{b}
\Bput{\Id}
\ncline{d}{c}
\Bput{\Id}
\ncline{f}{e}
\Bput{\Id}
\ncline{h}{g}
\Bput{\Id}
\psset{linestyle=dotted,arrows=-}
\ncline{j}{c}
\ncline{d}{e}
\ncline{f}{g}
\end{equation}

\vspace*{1.8cm}

\noindent
where $W'=W/\text{im} ([\YY,\XX]-\Id)$ and $\ZZ:W\r W'$ is the quotient map.

Let $E$ be the right $A$-module which is the restriction of
$\Escr$. Since $\Escr$ has a left $B$ structure it follows that $E$
is a $B-A$-bimodule. As a $\Mod(\Delta^{\opp})-A$ object it is given
by 

\vspace*{1.1cm}

\[
\psset{arrows=->,nodesep=0.1cm,labelsep=0.1cm,border=1mm}
\arraycolsep=0.4cm
\begin{array}{cccccccccccc}
\rnode{a}{A}  & \rnode{b}{A} & \rnode{j}{A} & &\rnode{c}{A} & \rnode{d}{A} & &\rnode{e}{A} & \rnode{f}{A} & &
\rnode{g}{A}
& \rnode{h}{A}
\end{array}
\everypsbox{\scriptstyle}
\ncarc[arcangle=45]{a}{c}
\Aput{x\cdot}
\ncarc[arcangle=45]{b}{d}
\Aput{x\cdot}
\ncarc[arcangle=45]{e}{h}
\Aput{x\cdot}
\ncarc[arcangle=-45]{a}{e}
\bput{0}(0.4){y\cdot}
\ncarc[arcangle=-45]{b}{f}
\bput{0}(0.55){y\cdot}
\ncarc[arcangle=-45]{c}{h}
\Bput{y\cdot}
\ncline{a}{b}
\Aput{\Id}
\ncline{b}{j}
\Aput{\Id}
\ncline{c}{d}
\Aput{\Id}
\ncline{e}{f}
\Aput{\Id}
\ncline{g}{h}
\Aput{\Id}
\psset{linestyle=dotted,arrows=-}
\ncline{c}{j}
\ncline{e}{d}
\ncline{g}{f}
\]

\vspace*{1.8cm}

It is clear from the above discussion that the inverse to
\eqref{app1} is given by 
\begin{equation}
\label{leftinverse}
(W,\XX,\YY)\mapsto \Tor_1^B(M,E)
\end{equation}
Let $\proj(B)$ be the category of finitely generated projective right
$B$-modules and let $\Sigma$ be the collection of maps 
$Z_1, \ldots, Z_{|w|-1}$. Let $B_{\Sigma}$ be the universal localization of
$B$ at $\Sigma$. 

Let us recall how $B_\Sigma$ is constructed \cite{schofield1}. We adjoin
 the inverses of the maps in $\Sigma$ to $\proj(B)$. Denote the resulting  category
by $\Sigma^{-1}\proj(B)$. Then
$B_\Sigma=\End_{\Sigma^{-1}\proj(B)}(B)$. Now in our case
$\proj(B)\cong \Delta$ and under this equivalence $B$ corresponds to
the sum of the vertices. It is also clear that $\Sigma^{-1}\proj(B)$
is equivalent to  $\Sigma^{-1}\Delta$ which is obtained from
$\Delta$ by adjoining inverses to the arrows in $\Sigma$. Then
$B_\Sigma$ is the endomorphism ring of the sum of the vertices in
$\Sigma^{-1}\Delta$, i.e. the path algebra. Thus we obtain 
$B_\Sigma= \CC(\Sigma^{-1}\Delta)$.  

Now let $\Delta^0$ be the following quiver.

\vspace*{0.5cm}

\[
\psset{arrows=->,nodesep=0.1cm,labelsep=0.1cm}
\arraycolsep=1cm
\begin{array}{cc}
\rnode{a}{0} & \rnode{b}{1}
\end{array}
\everypsbox{\scriptstyle}
\ncline[arrows=->]{a}{b}
\Aput{Z_0}
\psset{arrows=-}
\ncarc[arcangleA=60, arcangleB=60, ncurv=5]{b}{b}
\Bput{X_1}
\ncarc[arcangleA=-120, arcangleB=-120, ncurv=5]{b}{b}
\Bput{Y_1}
\]

\vspace*{0.7cm}

\noindent
with relation $(X_1Y_1-Y_1X_1-\Id)Z_0=0$. The obvious functor $\Sigma^{-1}\Delta\r \Delta^0$ which sends the arrows in
$\Sigma$ to the identity on the vertex $1$ is an equivalence of
categories. So we  find
$\Mod(\Sigma^{-1}\Delta)\cong\Mod(\Delta^0)$.
Below we put $B^0=\CC(\Delta^0)$.

Now we return to \eqref{leftinverse}. It is clear from the quiver
description \eqref{Mdescription} that $M$ may be
viewed (necessarily in unique way) as an object in
$\Mod(\Sigma^{-1}\Delta)$. Hence $M$ is a right $B_\Sigma$-module.  In
a similar way it follows that $E$ is a  $B_\Sigma-A$-bimodule. Then
according to \cite[Thm 4.8(c)]{schofield1} we have
$\Tor_1^B(M,E)=\Tor_1^{B_\Sigma}(M,E)$. 

Under the equivalence $\Mod(\Delta^0)\cong \Mod(\Sigma^{-1}\Delta)$
the right $B_\Sigma$-module $M$ corresponds to $M^0$ which is given by 

\vspace*{0.5cm}

\[
\psset{arrows=->,nodesep=0.1cm,labelsep=0.1cm}
\arraycolsep=1cm
\begin{array}{ccc}
\rnode{x}{M^0:}&\rnode{a}{W'} & \rnode{b}{W}
\end{array}
\everypsbox{\scriptstyle}
\ncline[arrows=->]{b}{a}
\Aput{\ZZ}
\psset{arrows=-}
\ncarc[arcangleA=60, arcangleB=60, ncurv=5]{b}{b}
\Bput{\XX}
\ncarc[arcangleA=-120, arcangleB=-120, ncurv=5]{b}{b}
\Bput{\YY}
\]

\vspace*{0.7cm}

Now $\Tor_1^{B_\Sigma}(-,E)$
is the first left derived functor of the functor $-\otimes_{B_\Sigma}
E:\Mod(\Sigma^{-1}\Delta)\r \Mod(A)$. Checking on projectives we see
that if we compose this functor with the equivalence
$\Mod(\Delta^0)\cong \Mod(\Sigma^{-1}\Delta)$ then it is given by
tensoring with the $B^0-A$-module $E^0$ which is defined as follows:

\vspace*{0.5cm}

\[
\psset{arrows=->,nodesep=0.1cm,labelsep=0.1cm}
\arraycolsep=1cm
\begin{array}{ccc}
\rnode{x}{E^0:}&\rnode{a}{A} & \rnode{b}{A}
\end{array}
\everypsbox{\scriptstyle}
\ncline[arrows=->]{a}{b}
\Aput{\Id}
\psset{arrows=-}
\ncarc[arcangleA=60, arcangleB=60, ncurv=5]{b}{b}
\Bput{x\cdot}
\ncarc[arcangleA=-120, arcangleB=-120, ncurv=5]{b}{b}
\Bput{y\cdot}
\]

\vspace*{0.7cm}

Thus we have shown that the inverse to \eqref{app1} is given by
the functor
\begin{equation}
\label{leftinverse1}
(W,\XX,\YY)\mapsto \Tor_1^{B^0}(M^0,E^0)
\end{equation}
It is clear that this inverse does not depend on $w$. This finishes
the proof of Theorem 1.4.

\medskip

We will now prove Theorem 1.3 by showing that \eqref{leftinverse1} is
compatible with the $\Aut(A)$-actions. Let us recall how these actions
are defined. $\Aut(A)$ is generated by the automorphisms
$\Psi_{n,\lambda}$ and $\Phi_{m,\mu}$ defined in the
introduction. As explained in  Section  7 there is an $\Aut(A)$ action on
$\Cscr_2$ which on the generators $\Psi_{n,\lambda}$ and
$\Phi_{m,\mu}$ is given by
$\Psi_{n,\lambda}(\XX,\YY)=(\XX,\YY-\lambda \XX^n)$ and
$\Phi_{m,\mu}(\XX,\YY)=(\XX-\mu \YY^m,\YY)$. 

We also define an action of $\Aut(A)$ on $\Delta^0$ (viewed as an
additive category) by 
\begin{align*}
\Psi_{n,\lambda}(X^1)&=X^1 & \Phi_{m,\mu}(X^1)&=X^1+\mu (X^1)^m\\
\Psi_{n,\lambda}(Y^1)&=Y^1+\lambda (X^1)^n & \Phi_{m,\mu}(Y^1)&=Y^1\\
\Psi_{n,\lambda}(Z^0)&=Z^0 & \Phi_{m,\mu}(Z^0)&=Z^0
\end{align*}
This is well defined because of the Remark in Section  7. Thus we
obtain an action of $\Aut(A)$ on $B^0$ in the obvious way. We obtain a
corresponding action on $\Cscr_1$ by putting  $\sigma(M^0)=M^0_{\sigma^{-1}}$
where $\sigma\in \Aut(A)$ and $M^0_{\sigma^{-1}}$ is the right
$B^0$-module which is equal to $M^0$ as a set but whose right
$B^0$-action is twisted by $\sigma^{-1}$. The action of $\Aut(A)$ on
$R$ is
defined similarly.

By checking on the
generators  $\Psi_{n,\lambda}$ and
$\Phi_{m,\mu}$ it is easy to see that the actions of
$\Aut(A)$ on $\Cscr_1$ and $\Cscr_2$ are compatible. 
Hence to prove our claim it is sufficient to prove that the functor $M^0\mapsto
\Tor_1^{B^0}(M^0,E)$ is compatible with the $\Aut(A)$-actions. To prove this
we need that 
\begin{equation}
\label{claim}
{}_\sigma E\cong E_{\sigma^{-1}}
\end{equation}
as $B^0-A$-bimodules, since if this is the case then $\Tor_1^{B^0}(M^0_{\sigma^{-}},E)\cong
\Tor_1^{B^0}(M^0,{}_\sigma E)\cong\Tor_1^{B^0}(M^0, 
E_{\sigma^{-1}})\cong\Tor_1^{B^0}(M^0, E)_{\sigma^{-1}}$. 

To prove \eqref{claim} we note that ${}_\sigma E$ is the
$\Mod(\Delta^{\opp})-A$ object given by the top quiver in the diagram
below.

\vspace*{1cm}

\[
\psset{arrows=->,nodesep=0.1cm,labelsep=0.1cm,border=1mm}
\arraycolsep=1cm
\begin{array}{cc}
\rnode{a}{A} & \rnode{b}{A}\\[2cm]
\rnode{a1}{A} & \rnode{b1}{A}
\end{array}
\ncline[linestyle=dotted]{a}{a1}
\Aput{\sigma^{-1}}
\ncline[linestyle=dotted]{b}{b1}
\Aput{\sigma^{-1}}
\ncline{a}{b}
\Aput{\Id}
\psset{arrows=-}
\ncarc[arcangleA=-30, arcangleB=150, ncurv=5]{b}{b}
\Bput{\sigma(x)\cdot}
\ncarc[arcangleA=60, arcangleB=60, ncurv=5]{b}{b}
\Bput{\sigma(y)\cdot}
\ncarc[arcangleA=-30, arcangleB=150, ncurv=5]{b1}{b1}
\Bput{x\cdot}
\ncarc[arcangleA=-120, arcangleB=-120, ncurv=5]{b1}{b1}
\Bput{y\cdot}
\psset{arrows=->}
\ncline{a1}{b1}
\Bput{\Id}
\]

\vspace*{1cm}

\noindent It is clear that the map given by the dotted arrows is left
$B^0$-linear and that it twists the right $A$-action by $\sigma^{-1}$.
This finishes the proof of Theorem 1.3.

\ifx\undefined\bysame
\newcommand{\bysame}{\leavevmode\hbox to3em{\hrulefill}\,}
\fi
\bibliographystyle{amsalpha}

\end{document}